\documentclass[11pt,a4paper]{article}
\usepackage{amsmath,amssymb,amsthm,mathrsfs}

\renewcommand{\mathcal}{\mathscr}

\renewcommand{\P}{\mathrm{P}}
\newcommand {\E}{{\mathrm E}}
\newcommand{\1}{{\bf 1}}

\newcommand{\R}{\mathbb{R}}
\newcommand{\e}{\epsilon}

\newtheorem{stat}{Statement}[section]
\newtheorem{prop}[stat]{Proposition}
\newtheorem{cor}[stat]{Corollary}
\newtheorem{thm}[stat]{Theorem}
\newtheorem{lem}[stat]{Lemma}

\theoremstyle{definition} \newtheorem{rem}[stat]{Remark}
\numberwithin{equation}{section}

\begin{document}
\title{\sc The fractional stochastic heat equation on the circle: Time regularity and potential theory}

\author{Eulalia Nualart$^1$ and Frederi Viens$^{2,3}$}

\date{}

\maketitle

\footnotetext[1]{Institut Galil\'ee, Universit\'e
Paris 13, 93430 Villetaneuse, France. 
\texttt{nualart@math.univ-paris13.fr} }

\footnotetext[2]{Department of Statistics, Purdue 
	University, 150 N. University St., West Lafayette, IN 47907-2067, USA. \texttt{viens@purdue.edu}}

\footnotetext[3]{The research of this author is partially supported by NSF grant num. 0606615-DMS}
 
\maketitle
\begin{abstract}  
We consider a system of $d$ linear stochastic
	heat equations driven by an additive infinite-dimensional fractional Brownian noise on the unit circle $S^1$. 
We obtain sharp results on the H\"older continuity in time of the paths of the solution $u=\{u(t\,,x)\}_{t \in 
	\mathbb{R}_+, x \in S^1}$. We then establish
	upper and lower bounds on hitting
	probabilities of $u$, in terms of respectively
	Hausdorff measure and Newtonian capacity.
	\end{abstract}

\vskip 1,5cm {\it \noindent AMS 2000 subject classifications:}
Primary: 60H15, 60J45; Secondary: 60G15, 60G17. \vskip 10pt


\vskip 6cm \pagebreak

\section{Introduction and main results}

We consider a system of $d$ stochastic heat equations on the unit circle
driven by an infinite-dimensional fractional Brownian motion $B^{H}$ with
Hurst parameter $H\in (0,1)$. That is, 
\begin{equation}
\frac{\partial u_{i}}{\partial t}(t,x)=\Delta _{x}u_{i}(t,x)+\frac{\partial
B_{i}^{H}}{\partial t}(t,x),\;\;t>0,\;\;x\in S^{1},  \label{equa1}
\end{equation}
with initial condition $u_{i}(0,x)=0$, for all $i=1,...,d$. Here $\Delta
_{x} $ is the Laplacian on $S^{1}$ and $B^{H}$ a centered Gaussian field on $
\mathbb{R}_{+}\times S^{1}$ defined, for all $x,y\in S^{1}$ and $s,t\geq 0$,
by its covariance structure
\begin{equation*}
\E\left[ B_{i}^{H}\left( t,x\right) B_{j}^{H}\left( s,y\right) 
\right] =2^{-1}\left( t^{H}+s^{H}-\left\vert t-s\right\vert ^{2H}\right)
Q\left( x,y\right) \delta _{i,j},
\end{equation*}
where $Q$ is an arbitrary covariance function on $S^{1}$ and $\delta _{i,j}$
is the Kronecker symbol. To simplify our study, we assume that $B^{H}$ is
spatially homogeneous and separable in space; therefore $Q\left( x,y\right) $
depends only on the difference $x-y$ , and we denote it abusively $Q(x-y)$.

Note that because $Q$ is positive definite, there exists a sequence of
non-negative real numbers $\left\{ q_{n}\right\} _{n\in \mathbb{N}}$ such
that 
\begin{equation*}
Q\left( x-y\right) =\sum_{n\in \mathbb{N}}q_{n}\cos \left( n(x-y)\right) .
\end{equation*}
This expression may be only formal for certain choices of the sequence $
\left\{ q_{n}\right\}_{n}$, as these pointwise values may explode, but this
Fourier representation is always relevant if one allows $Q$ to be a Schwartz
distribution. Examples will be given below where $Q\left( 0\right) $ is
infinite while all other values are finite (Riesz kernel case); another,
also with $Q\left( 0\right) =\infty $, will show that $Q$ may not be equal
to its Fourier series at any point (fractional noise case for small Hurst
parameter), but still allows a solution to (\ref{equa1}). Any case with $
Q\left( 0\right) =\infty $ denotes a distribution-valued noise $B^{H}$ in
space, for which the notation $B^H\left( t,x\right) $ is only formal in the
parameter $x$.

The existence and uniqueness of the solution of (\ref{equa1}) was
established in \cite{Tindel:03}. The \textquotedblleft
mild\textquotedblright\ or \textquotedblleft evolution\textquotedblright\
solution of the stochastic integral formulation of equation (\ref{equa1}) is
given by the evolution convolution
\begin{equation}
u_{i}\left( t,x\right) =\sum_{n=0}^{\infty }\sqrt{q_{n}}\biggl(\cos \left(
nx\right) \int_{0}^{t}e^{-n^{2}\left( t-s\right) }\beta _{i,n}^{H}\left(
ds\right) +\sin \left( nx\right) \int_{0}^{t}e^{-n^{2}\left( t-s\right)
}\beta _{i,n}^{\prime H}\left( ds\right) \biggr),  \label{equa2}
\end{equation}
where the sequences $\{\beta _{i,n}\}_{n\in \mathbb{N}}$ and $\{\beta
_{i,n}^{\prime }\}_{n\in \mathbb{N}}$, $i\in \{1,...,d\}$, are independent
and each formed of independent one-dimensional standard fractional Brownian
motions. \cite{Tindel:03} showed when such a solution exist, and more
specifically, that the necessary and sufficient condition for existence of 
(\ref{equa2}) in $L^{2}\left( \Omega \times \lbrack 0,T]\times S^{1}\right) $
(cf.~\cite[Corollary 1]{Tindel:03}) is 
\begin{equation*}
\sum_{n=1}^{\infty }q_{n}n^{-4H}<\infty .
\end{equation*}
The particular case where $B^{H}$ is white in space, that is, $q_{n}=1$ for
all $n$, was already studied in \cite{Duncan:02} where the solution exists
if and only if $H>\frac{1}{4}$.

The aim of this paper is to develop a potential theory for the solution to
the system of equations (\ref{equa1}). In particular, given $A\subset 
\mathbb{R}^{d}$, we want to determine whether the process $
\{u(t\,,x),\,t\geq 0,\,x\in S^{1}\}$ visits, or hits, $A$ with positive
probability.

Potential theory for the linear and non-linear stochastic heat equation
driven by a space time white noise was developped in \cite{Dalang:05} and 
\cite{Dalang2:05}. The aim of this paper is to obtain upper and lower bounds
on hitting probabilities for the solution of (\ref{equa1}). For this,
following the approach developped in \cite{Dalang:05} , a careful analysis
of the moments of the increments of the process $u(t,x)$ is needed. In
particular, this will lead us to solve an open question which is the 
H\"{o}lder continuity in time of the solution of (\ref{equa1}) when $H<\frac{1}{2}$. The H\"older continuity in space for the solution of (\ref{equa1}) was
studied in \cite{Tindel:04} and the H\"older continuity in time when $
H\geq \frac{1}{2}$ is due to \cite{Sarol:06}. These are generalizations of earlier work 
done for the stochastic heat equation with time-white noise potential: \cite{Sanz:00}, \cite{Sanz:02}.

Let us first state, in some detail, the path continuity results we obtain
for the solution of the fractional heat equation on the circle (\ref{equa1}), as these are a valuable immediate consequence of our work. Assume that
for all $n$ large enough 
\begin{equation}
cn^{4H-2\alpha -1}\leq q_{n}\leq Cn^{4H-2\alpha -1},  \label{hq}
\end{equation}
for some positive constants $c$ and $C$ and $\alpha \in (0,1]$ with $\alpha \neq 2H$. Our basic quantitative result is the following bounds on the
variance of the increments of the solution: for $t_{0},T>0$, for some
positive constants $c,C,c_{t_{0}},C_{t_{0}}$, for all $x,y\in S^{1}$, and
all \thinspace $s,t\in \lbrack t_{0},T]$,
\begin{eqnarray*}
c_{t_{0}}\left\vert x-y\right\vert ^{2\alpha } &\leq &\E\left[
\Vert u\left( t,x\right) -u\left( t,y\right) \Vert ^{2}\right] \leq
C_{t_{0}}\left\vert x-y\right\vert ^{2\alpha } \\
c\left\vert t-s\right\vert ^{\alpha \wedge (2H)} &\leq &\E\left[
\Vert u\left( t,x\right) -u\left( s,x\right) \Vert ^{2}\right] \leq
C\left\vert t-s\right\vert ^{\alpha \wedge (2H)}.
\end{eqnarray*}

We then immediately get that $u$ is $\beta $-H\"{o}lder continuous in space
for any $\beta \in (0,\alpha )$ and is $\beta $-H\"{o}lder continuous in
time for any $\beta \in (0,\frac{\alpha }{2}\wedge H)$, but not for $\beta$
equal to the upper values of these intervals. All these results are true for
any $H\in (0,1)$. Moreover, these results are sharp for our additive
stochastic heat equation (\ref{equa1}): up to non-random constants, exact
moduli of continuity can be found (see the last bullet point below).

Let us consider some examples:
\begin{itemize}

\item In the case where $B^{H}$ is \textquotedblleft white
noise\textquotedblright\ in space, then $u$ exists if and only if $H>1/4$;
moreover $u$ is $\beta $-H\"{o}lder continuous in space for any $\beta \in
(0,2H-\frac{1}{2})$ and $\beta $-H\"{o}lder continuous in time for any $
\beta \in (0,H-\frac{1}{4})$. This follows from the above continuity results
because the white noise case is the case $q_{n}\equiv 1$: the appellation
\textquotedblleft white\textquotedblright\ reflects the fact that all
spatial Fourier frequencies are equally represented.

\item In the case where $B^{H}$ is white in time and has a covariance
function in space given by the Riesz kernel, that is, $Q(x-y)=|x-y|^{-\gamma
}$, $0<\gamma <1$, we can prove that $q_{n}$ is commensurate with $n^{\gamma
-1}$. More specifically, we can show that $q_{n}=n^{\gamma -1}c\left(
n\right) $ where $c\left( n\right) $ is function bounded bewteen two
positive constants, because it can be written as the partial sum of an
alternating series with decreasing general term and positive initial term
(see Appendix \ref{a11}). Therefore, the solution of (\ref{equa1}) exists if and
only if $H>\frac{\gamma }{4}$ and $u$ is $\beta $-H\"{o}lder continuous in
space for any $\beta \in (0,2H-\frac{\gamma }{2})$ and $\beta $-H\"{o}lder
continuous in time for any $\beta \in (0,H-\frac{\gamma }{4})$.

\item In the case where $B^{H}$ behaves like fractional Brownian noise both
in time and space with common Hurst parameter $H$, then the solution of (\ref
{equa1}) exists if and only if $H>\frac{1}{3}$. Indeed, this case can be
obtained by assuming that $q_{n}=n^{1-2H}$. When $H>1/2$, if one prefers to
work starting from the spatial covariance function $Q$, one may stipulate
that $B^{H}$ is has a Riesz kernel covariance, i.e. $Q\left( x-y\right)
=\left\vert x-y\right\vert ^{2H-2}=\left\vert x-y\right\vert ^{-\gamma }$
with $\gamma =2\left( 1-H\right) \in (0,1)$, in which case one is in the
situation of the last example, with $q_{n}=c\left( n\right) n^{1-2H}$. On
the other hand, if $H\leq 1/2$, no Riesz-kernel interpretation is possible
with $q_{n}=n^{1-2H}$. Appendix \ref{a2} contains another interpretation in this
case.

This interpretation, which uses a differentiation construction, also allows
a justification, for all $H\in (0,1)$, of why we use the appellation
\textquotedblleft fractional Brownian noise\textquotedblright\ in the case $
q_{n}=n^{1-2H}$. In all cases, i.e. for all $H\in (1/3,1)$, $u$ is $\beta $-H
\"{o}lder continuous in space for any $\beta \in (0,3H-1)$ and is $\beta $-H
\"{o}lder continuous in time for any $\beta \in (0,\frac{3H-1}{2})$.

\item Similarly to the previous example, but more generally, to obtain a $
B^{H}$ that behaves like a fractional Brownian noise with parameter $H$ in
time and $K$ in space, we can set $q_{n}=n^{1-2K}$ (using the same
justification as in the Appendix relative to the previous example). This is
equivalent to $\alpha =2H+K-1$. We then get existence of a solution if and
only if $2H+K>1$, and the solution is then $\beta $-H\"{o}lder continuous in
space for any $\beta \in (0,2H+K-1)$ and is $\beta $-H\"{o}lder continuous
in time for any $\beta \in (0,\frac{2H+K-1}{2})$.

\item From Gaussian regularity results such as Dudley's entropy upper bound (see \cite{Khoshnevisan:02}), 
we can state that
if the upper bound in (\ref{hq}) holds, then the modulus of continuity
random variable
\begin{equation*}
\sup_{x,y\in S^{1};s,t\in \lbrack t_{0},T]}\left( \frac{\left\Vert u\left(
t,x\right) -u\left( t,y\right) \right\Vert }{\left\vert x-y\right\vert
^{\alpha }\log ^{1/2}\left( 1+1/\left\vert x-y\right\vert \right) }+\frac{
\left\Vert u\left( t,x\right) -u\left( s,x\right) \right\Vert }{\left\vert
t-s\right\vert ^{(\alpha /2)\wedge H}\log ^{1/2}\left( 1+1/\left\vert
t-s\right\vert \right) }\right)
\end{equation*}
is finite almost surely. Moreover, a (near) converse also holds: if the
above random variable (with logarithmic terms moved to the numerators) is
finite, then the upper bound in (\ref{hq}) holds for some constant $C<\infty$ 
(see \cite[Corollary 1]{Tindel:04}). 
\end{itemize}

We now state the results of potential theory that we will prove in this
paper. For this, let us first introduce some notation. For all Borel sets $
F\subseteq\mathbb{R}^d$ we define $\mathcal{P}(F)$ to be the set of all
probability measures with compact support in $F$. For all $\mu\in\mathcal{P}(
\mathbb{R}^d)$, we let $I_\beta(\mu)$ denote the $\beta$-dimensional energy
of $\mu$; that is, 
\begin{equation*}
I_\beta(\mu) := \iint \mathrm{K}_\beta(\|x-y\|)\, \mu(dx)\,\mu(dy).
\end{equation*}
Here and throughout, 
\begin{equation}  \label{k}
\mathrm{K}_\beta(r) := 
\begin{cases}
r^{-\beta} & \text{if $\beta >0$}, \\ 
\log ( N_0/r ) & \text{if $\beta =0$}, \\ 
1 & \text{if $\beta<0$},
\end{cases}
\end{equation}
where $N_0$ is a constant whose value will be specified later in the proof
of Lemma \ref{prel}.

For all $\beta\in\mathbb{R}$ and Borel sets $F\subset\mathbb{R}^d$, $\text{
Cap}_\beta(F)$ denotes the $\beta$-dimensional capacity of $F$; that is, 
\begin{equation*}
\text{Cap}_\beta(F) := \left[ \inf_{\mu\in\mathcal{P}(F)} I_\beta(\mu) 
\right]^{-1},
\end{equation*}
where $1/\infty:=0$.

Given $\beta\geq 0$, the $\beta$-dimensional Hausdorff measure of $F$ is defined by 
\begin{equation*}  \label{eq:HausdorfMeasure}
{\mathcal{H}}_ \beta (F)= \lim_{\epsilon \rightarrow 0^+} \inf \left\{
\sum_{i=1}^{\infty} (2r_i)^ \beta : F \subseteq \bigcup_{i=1}^{\infty}
B(x_i\,, r_i), \ \sup_{i\ge 1} r_i \leq \epsilon \right\},
\end{equation*}
where $B(x\,,r)$ denotes the open (Euclidean) ball of radius $r>0$ centered
at $x\in \mathbb{R}^d$. When $\beta <0$, we define $\mathcal{H}_\beta (F)$
to be infinite.

Let $u(S)$ denote the range of $S$ under the random map $r\mapsto u(r)$,
where $S$ is some Borel-measurable subset of $\mathbb{R}_{+}\times S^{1}$.
\begin{thm}\label{t1} 
Assume hypothesis \textnormal{(\ref{hq})}.
Let $I\subset (0,T]$ and $J\subset \lbrack 0,2\pi )\equiv S^{1}$ be two
fixed non-trivial compact intervals. Then for all $T>0$ and $M>0$, there exists a finite constant $c_{H}>0$ depending on $
H,M,I$ and $J$ such that for all compact sets $A\subseteq \lbrack -M,M]^{d}$, 
\begin{equation*}
c_{H}^{-1}\,\textnormal{Cap}_{d-\beta }(A)\leq \P\{u(I\times J)\cap
A\neq \emptyset \}\leq c_{H}\,\mathcal{H}_{d-\beta }(A).
\end{equation*}
where $\beta :=\frac{1}{\alpha }+(\frac{2}{\alpha }\vee \frac{1}{H})$.
\end{thm}

\begin{rem}
\begin{itemize}
\item[\textnormal{(a)}] When $B^{H}$ is white in time and space, that is, $H=
\frac{1}{2}$ and $q_{n}=1$ for all $n$, Theorem \ref{t1} gives the same
hitting probabilities estimates obtained in \cite[Theorem 4.6.]{{Dalang:05}}.

\item[(b)] Because of the inequalities between capacity and Hausdorff
measure, the right-hand side of Theorem \ref{t1} can be replaced by $c\,
\textnormal{Cap}_{d-\beta -\eta }(A)$ for all $\eta >0$ (cf. \cite[p. 133]{Kahane:85}).
\end{itemize}
\end{rem}

We say that a Borel set $A\subseteq \mathbb{R}^{d}$ is called polar for $u$
if $\P\{u(T)\cap A\neq \varnothing \}=0$; otherwise, $A$ is called
nonpolar.

The following results are consequences of Theorem \ref{t1}.
\begin{cor}\label{c1} 
Assume hypothesis \textnormal{(\ref{hq})} and let $\beta:=\frac{1
}{\alpha}+(\frac{2}{\alpha} \vee \frac1H)$.
\begin{itemize}
\item[\textnormal{(a)}] A (nonrandom) Borel set $A\subset\mathbb{R}^d$ is
nonpolar for u if it has positive $d-\beta$-dimensional capacity. On the
other hand, if $A$ has zero $d-\beta$-dimensional Hausdorff measure, then $A$
is polar for $u$.

\item[\textnormal{(b)}] Singletons are polar for $u$ if $d > \beta$ and are
nonpolar when $d<\beta$. The case $d=\beta$ is open.

\item[\textnormal{(c)}] If $d\geq \beta$, then 
\begin{equation*}
\dim_{_{\mathrm{H}}} ({u}( \mathbb{R}_+ \times\, S^1)) = \beta , \; \; 
\text{a.s.}
\end{equation*}
\end{itemize}
\end{cor}

Let us consider the same examples as we had for the regularity statements.
\begin{itemize}
\item In the case where $B^{H}$ is white in space, then $\alpha =2H-\frac{1}{
2}$ and $\beta =\frac{6}{4H-1}$.

\item In the case where $B^{H}$ is white in time and has a covariance
function in space given by the Riesz kernel, that is, $Q(x-y)=|x-y|^{-\gamma
}$, $0<\gamma <1$, then $\alpha =2H-\frac{\gamma }{2}$ and $\beta =\frac{6}{
4H-\gamma }$.

\item In the case where $B^{H}$ is the fractional Brownian noise with Hurst
parameter $H>1/3$ in time and space, then $\alpha =3H-1$ and $\beta =\frac{3
}{3H-1}$.

\item In the case where $B^{H}$ is the fractional Brownian noise with Hurst
parameter $H$ in time and $K$ in space, and $2H+K>1$, then $\alpha =2H+K-1$
and $\beta =\frac{3}{2H+K-1}$.
\end{itemize}

This paper is organized as follows. In Section 2 we prove the path continuity results of $
u$ stated in the Introduction using fractional stochastic calculus.
In Section 3 we obtain an upper bound of Gaussian type for the bivariate
density of $u$ that will be needed for the proof of Theorem \ref{t1}. Finally, Section 4 is devoted to the proofs of Theorem \ref
{t1} and Corollary \ref{c1}. 

In all the paper, $c_H, C_H$ will denote universal constants depending on $H$
whose value may change from line to line.

\section{Regularity of the solution}

We consider the two canonical metrics of $u$ in the space and time
parameter, respectively, defined by 
\begin{equation*}
\begin{split}
\delta^2_t(x,y)&:= {\mathrm{E}} [\Vert u(t,x)-u(t,y) \Vert^2], \\
\delta^2_x(s,t)&:= {\mathrm{E}} [\Vert u(t,x)-u(s,x) \Vert^2],
\end{split}
\end{equation*}
for all $x,y \in S^1$ and $s,t \in \mathbb{R}_+$.

The aim of this section is to obtain upper and lower bounds bounds in terms
of the differences $|x-y|$ and $|t-s|$ for the two canonical metrics above.
These imply, in particular, the H\"{o}lder regularity of $u$ that we have
described in detail in the introduction. We begin by introducing some
elements of fractional stochastic calculus.

\subsection{Elements of fractional stochastic calculus}

In this section, we recall, following \cite{Nualart:06}, some elements on
stochastic integration with respect to one-dimensional fractional Brownian
motion needed for the analysis of the regularity of $u$ in time.

Fix $T>0$. Let $B^H=(B^{H}(t),t\in \lbrack 0,T])$ be a one-dimensional
fractional Brownian motion with Hurst parameter $H\in (0,1)$. That is, $B^H$
is a centered Gaussian process with covariance function given by 
\begin{equation*}
R(t,s)={\mathrm{E}}[B^{H}(t)B^{H}(s)]=2^{-1}\left( t^{H}+s^{H}-\left\vert
t-s\right\vert ^{2H}\right) .
\end{equation*}
Note that for $H=\frac{1}{2}$, $B^{H}$ is a standard Brownian motion.
Moreover, $B^{H}$ has the integral respresentation 
\begin{equation*}
B^{H}(t)=\int_{0}^{t}K^{H}(t,s)W(ds),
\end{equation*}
where $W=(W(t),t\in \lbrack 0,T])$ is a Wiener process and $K^{H}(t,s)$ is
the kernel defined as 
\begin{equation}
K^{H}(t,s)=c_{H}\biggl(\frac{t}{s}\biggr)^{H-\frac{1}{2}}(t-s)^{H-\frac{1}{2}
}+s^{\frac{1}{2}-H}F\biggl(\frac{t}{s}\biggr),  \label{kernel}
\end{equation}
where $c_{H}$ is a positive constant and 
\begin{equation*}
F(z)=c_{H}\biggl(\frac{1}{2}-H\biggr)\int_{0}^{z-1}r^{H-\frac{3}{2}}\biggl(
1-(1+r)^{H-\frac{1}{2}}\biggr)dr.
\end{equation*}
From (\ref{kernel}) we get 
\begin{equation}
\frac{\partial K^{H}}{\partial t}(t,s)=c_{H}\biggl(H-\frac{1}{2}\biggr)
(t-s)^{H-\frac{3}{2}}\biggl(\frac{s}{t}\biggr)^{\frac{1}{2}-H}.
\label{dkernel}
\end{equation}
It is important to note that $\frac{\partial K^{H}}{\partial t}$ is positive
if $H>1/2$, but is negative when $H<1/2$. This negativity causes problems
when evaluating the time-canonical metric's lower bound.

We denote by $\mathcal{E}$ the set of step functions on $[0,T]$. Let $
\mathcal{H}$ be the Hilbert space defined as the closure of $\mathcal{E}$
with respect to the scalar product 
\begin{equation*}
\langle \mathbf{1}_{[0,t]},\mathbf{1}_{[0,s]}\rangle _{\mathcal{H}}=R(t,s).
\end{equation*}
The mapping $\mathbf{1}_{[0,t]}\mapsto B_{t}^{H}$ can be extended to an
isometry between $\mathcal{H}$ and the Gaussian space $\mathcal{H}_{1}$
associated with $B^{H}$. Then $\{B^{H}(\phi ),\phi \in \mathcal{H}\}$ is an
isonormal Gaussian process associated with the Hilbert space $\mathcal{H}$.
For every element $\phi \in \mathcal{H}$, $B^{H}(\phi )$ is called the
Wiener integral if $\phi $ with respect to $B^{H}$ and is denoted 
\begin{equation*}
\int_{0}^{T}\phi (s)B^{H}(ds).
\end{equation*}

For every $s<t$, consider the linear operator $K^{\ast }$ from $\mathcal{E}$
to $L^{2}([0,T])$ defined by 
\begin{equation*}
K_{t}^{\ast }\phi (s)=K^{H}(t,s)\phi (s)+\int_{s}^{t}(\phi (u)-\phi (s))
\frac{\partial K^{H}}{\partial u}(u,s)\,du.
\end{equation*}
When $H>\frac{1}{2}$, since $K^H\left( t,t\right) =0$, this operator has the
simpler expression 
\begin{equation*}
K_{t}^{\ast }\phi (s)=\int_{s}^{t}\phi (u)\frac{\partial K^{H}}{\partial u}
(u,s)\,du.
\end{equation*}
The operator $K^{\ast }$ is an isometry between $\mathcal{E}$ and $
L^{2}([0,T])$ that can be extended to the Hilbert space $\mathcal{H}$. As a
consequence, we have the following relationship between the Wiener integral
with respect to the fractional Brownian motion $B^{H}$ and the Wiener
integral with respect to the Wiener process $W$: 
\begin{equation*}
\int_{0}^{t}\phi (s)B^{H}(ds)=\int_{0}^{t}K_{t}^{\ast }\phi (s)W(ds),
\end{equation*}
which holds for every $\phi \in \mathcal{H}$, which is true if and only if $
K_{t}^{\ast }\phi \in L^{2}([0,T])$.

Recall also that when $H>\frac12$,
\begin{equation} \label{isometry}
\E \biggl[\int_{0}^{t}\phi (s)B^{H}(ds) \int_{0}^{t}\psi (s)B^{H}(ds)\biggr]=
H(2H-1) \int_{0}^{t} ds \int_{0}^{t} du \, \phi(s) \psi (u) \vert s-u\vert^{2H-2}.
\end{equation}

\subsection{Space regularity}

The next lemma gives a precise connection between a generic condition of the
type (\ref{hq}) and the Fourier expansion of a canonical metric for a
homogeneous Gaussian field on the circle.
\begin{lem}
\label{fBreg} Let $Y$ be a homogeneous, centered and separable Gaussian
field on $S^{1}$ with canonical metric $\delta \left( x,y\right) =\delta
\left( x-y\right) $ for some univariate function $\delta $. Then, there
exists a sequence of non-negative real numbers $\left\{ r_{n}\right\} _{n\in 
\mathbb{N}}$ such that for any $r\in S^{1}$, 
\begin{equation}
\delta ^{2}\left( r\right) =2\sum_{n=1}^{\infty }r_{n}\left( 1-\cos
nr\right) .  \label{deltafourier}
\end{equation}
Moreover, if there exist constants $c$ and $C$ positive, and $\alpha \in
(0,1]$, such that for all $n$ large enough, 
\begin{equation}
cn^{-2\alpha -1}\leq r_{n}\leq Cn^{-2\alpha -1},  \label{rnbounds}
\end{equation}
then for all $r$ close enough to $0$, 
\begin{equation}
\sqrt{k_{\alpha }c}r^{\alpha }\leq \delta \left( r\right) \leq \sqrt{
K_{\alpha }C}r^{\alpha },  \label{fBregbounds}
\end{equation}
where $k_{\alpha }$ and $K_{\alpha }$ are constants depending only on $
\alpha $. More specifically, the upper bound (resp. lower bound) in 
\textnormal{(\ref{rnbounds})} implies the upper bound (resp. lower bound) in 
\textnormal{(\ref{fBregbounds})}.
\end{lem}

\begin{proof}
We start proving (\ref{deltafourier}). Let $C(x,y)$ denote the covariance
function of $Y$, that is, for any $x,y\in S^{1}$, 
\begin{equation*}
{\mathrm{E}}[Y(x)Y(y)]=C(x,y),
\end{equation*}
where $C$ depends only on the diference $x-y$. Because $C$ is positive
definite, it holds that there exists a sequence of non-negative real numbers 
$\left\{ r_{n}\right\} _{n\in \mathbb{N}}$ such that 
\begin{equation*}
C(x,y)=\sum_{n\in \mathbb{N}}r_{n}\cos \left( n(x-y)\right) .
\end{equation*}
Hence, for any $r\in S^{1}$, 
\begin{equation*}
\delta ^{2}(r)={\mathrm{E}}[(Y(0)-Y(r))^{2}]=2\sum_{n=1}^{\infty
}r_{n}\left( 1-\cos nr\right) .
\end{equation*}
This proves (\ref{deltafourier}).

We now prove the second statement of the lemma. We begin proving the upper
bound statement. Assuming that the upper bound of (\ref{rnbounds}) holds for
all $n>n_{0}\geq 1$, we restrict $r$ accordingly: we assume $n_{0}\leq
\lbrack 1/r]$, that is, $r\leq 1/n_{0}$. In this case, we immediately get $
r^{2}\leq r^{2\alpha }$. We write 
\begin{align*}
2^{-1}\delta ^{2}\left( r\right) & =\sum_{n=1}^{n_{0}-1}r_{n}\left( 1-\cos
nr\right) +\sum_{n=n_{0}}^{[1/r]}r_{n}\left( 1-\cos nr\right)
+\sum_{n=[1/r]+1}^{\infty }r_{n}\left( 1-\cos nr\right) \\
& \leq \max_{n\leq n_{0}}\left\{ r_{n}\right\} \sum_{n=1}^{n_{0}-1}\left(
nr\right) ^{2}+\sum_{n=1}^{[1/r]}Cn^{-2\alpha -1}\left( nr\right)
^{2}+2\sum_{n=[1/r]+1}^{\infty }Cn^{-2\alpha -1} \\
& \leq n_{0}^{2}\max_{n\leq n_{0}}\left\{ r_{n}\right\}
r^{2}+Cr^{2}\sum_{n=1}^{[1/r]}n^{-2\alpha +1}+2\sum_{n=[1/r]+1}^{\infty
}Cn^{-2\alpha -1} \\
& \leq r^{2-2\alpha }n_{0}^{2}\max_{n\leq n_{0}}\left\{ r_{n}\right\}
r^{2\alpha }+CC_{\alpha }r^{2}\left( 1/r\right) ^{-2\alpha +2}+2CC_{\alpha
}^{\prime }\left( 1/r\right) ^{-2\alpha } \\
& \leq 2C\left( C_{\alpha }+2C_{\alpha }^{\prime }\right) r^{2\alpha },
\end{align*}
provided $r\leq r_{1}:=\min \left\{ 1/n_{0};C\left( C_{\alpha }+2C_{\alpha
}^{\prime }\right) \left[ n_{0}^{2}\max_{n\leq n_{0}}\left\{ r_{n}\right\} 
\right] ^{1/\left( 2-2\alpha \right) }\right\} $, where $C_{\alpha }$ and $
C_{\alpha }^{\prime }$ are constant depending only on $\alpha $. It is
elementary to check that $C_{\alpha }^{\prime }$ can be taken as $1/\left(
2\alpha \right) $. If $\alpha \in (0,1/2)$, then one checks that $C_{\alpha
} $ can be taken as $1$; while if $\alpha \in \lbrack 1/2,1]$, and we assume
moreover that $r<r_{2}:=\left( 1-2\alpha \right) ^{-1/\left( 2\alpha \right)
}$, then $C_{\alpha }$ can be taken as $\alpha ^{-1}$. In other words, when $
\alpha <1/2$, we obtain the upper bound of (\ref{fBregbounds}) for all $
r\leq r_{1}$, with $K_{\alpha }=4\left( \alpha ^{-1}+1\right) $, while when $
\alpha \in \lbrack 1/2,1]$, we obtain the upper bound of (\ref{fBregbounds})
for all $r\leq \min \left\{ r_{1};r_{2}\right\} $ with $K_{\alpha }=8\alpha
^{-1}$. In fact, the formula $K_{\alpha }=8\alpha ^{-1}$ can be used for
both cases.

In order to prove the lower bound on $\delta \left( r\right) $, we write
instead, still assuming $r\leq 1/n_{0}$, that 
\begin{align*}
2^{-1}\delta ^{2}\left( r\right) & =\sum_{n=1}^{\infty }r_{n}\left( 1-\cos
nr\right) \geq c\sum_{n=n_{0}}^{\infty }n^{-2\alpha -1}\left( 1-\cos
nr\right) \\
& \geq c\sum_{n=\left[ 1/r\right] +1}^{\left[ \pi /\left( 2r\right) \right]
}n^{-2\alpha -1}\left( 1-\cos nr\right) \geq c\left( 1-\cos 1\right) \sum_{n=
\left[ 1/r\right] +1}^{\left[ \pi /\left( 2r\right) \right] }n^{-2\alpha -1}\\
& \geq c\left( 1-\cos 1\right) \left( \frac{\pi }{2r}\right) ^{-2\alpha
-1}\left( \left[ \frac{\pi }{2r}\right] -1-\left[ \frac{1}{r}\right] \right)
\\
& \geq r^{2\alpha }c\left( 1-\cos 1\right) \left( \frac{\pi }{2}\right)
^{-2\alpha }\left( \frac{\pi }{2}-1-2r\right) .
\end{align*}

Note here that $1-\cos 1>0.459$ and $\pi /2-1>0.57$. It is now clear that
choosing $r\leq r_{0}:=\min \left\{ 0.035;1/n_{0}\right\} $, we get 
\begin{equation*}
\delta ^{2}\left( r\right) \geq r^{2\alpha }c\left( 1-\cos 1\right) \left(
\pi /2\right) ^{-2\alpha },
\end{equation*}
which proves the lower bound of (\ref{fBregbounds}) with $k_{\alpha }=\left(
1-\cos 1\right) \left( \pi /2\right) ^{-2\alpha }$ for all $r\leq r_{0}$.
The proof of the lemma is complete.
\end{proof}

This lemma can be applied immediately, to find sharp bounds on the spatial
canonical metric of $u$; the almost-sure continuity results also follow.
\begin{cor}
\label{corx} Let $H \in (0,1)$, $t_{0}>0$ and $t\in \lbrack t_{0},T]$ be fixed. Assume
hypothesis \textnormal{(\ref{hq})}. Then the canonical metric $\delta
_{t}\left( x-y\right) $ for $u\left( t,\cdot \right) $ satisfies, for all $r$
enough close to $0$, 
\begin{equation*}
\sqrt{k_{\alpha }cc\left( t_{0},T,H\right) }r^{\alpha }\leq \delta
_{t}\left( r\right) \leq \sqrt{K_{\alpha }CC\left( t_{0},T,H\right) }
r^{\alpha },
\end{equation*}
where $k_{\alpha }$ and $K_{\alpha }$ are constants depending only on $
\alpha $, $c\left( t_{0},T,H\right) $ and $c\left( t_{0},T,H\right) $ are
constant depending only on $t_{0},T$ and $H$ and $c,C$ are the constants in 
\textnormal{(\ref{hq})}.
In particular, $u(t,\cdot )$ is $\beta $-H\"{o}lder continuous for any $
\beta \in (0,\alpha )$. More specifically, up to a non-random constant, the
function $r\mapsto r^{\alpha }\log ^{1/2}\left( 1/r\right) $ is an almost
sure uniform modulus of continuity for $u\left( t,\cdot \right) $.
\end{cor}

\begin{proof}
Let $(\beta ^{H}\left( t\right) ,t\geq 0)$ be a one-dimensional fractional
Brownian motion. Let $t_{0}>0$ and $t\in \lbrack t_{0},T]$ be fixed. From
the proof of Theorems 2 and 3 of \cite{Tindel:03} we deduce that there
exists positive constants $c\left( t_{0},T,H\right) $ and $C\left(
t_{0},T,H\right) $ such that 
\begin{equation*}
c\left( t_{0},T,H\right) n^{-4H}\leq \mathrm{E}\left[ \left(
\int_{0}^{t}e^{-n^{2}\left( t-s\right) }\beta _{n}^{H}\left( ds\right)
\right) ^{2}\right] \leq C\left( t_{0},T,H\right) n^{-4H}.  
\end{equation*}
Thus, appealing to (\ref{equa2}), we find that for all $n$ sufficiently
large, 
\begin{equation*}
2c\left( t_{0},T,H\right) n^{-4H}q_{n}(1-\cos (nr))\leq \delta
_{t}^{2}(r)\leq 2C\left( t_{0},T,H\right) q_{n}n^{-4H}(1-\cos (nr))
\end{equation*}
Then hypothesis (\ref{hq}) and Lemma \ref{fBreg} conclude the first result
of the corollary.

The second statement of the corollary, which is a repeat of one of the
continuity results described in the introduction, is proved using the
arguments described therein as well. In fact, a simple application of 
Dudley's entropy upper bound theorem is sufficient (see \cite[Theorem 2.7.1]{Khoshnevisan:02}). We do not elaborate further on
this point.
\end{proof}

\subsection{Time regularity}

We now concentrate our efforts on finding sharp bound on the time-canonical
metric of $u$. The bounds we find for $H>1/2$ \ were essentially already
obtained in \cite{Sarol:06}, although the result and its proof was not
stated explicitly therein, an omission which we deal with here. When $H<1/2$, no results were known, either for upper or lower bounds: we perform these
calculations from scratch. This portion of our calculations is very
delicate. As in the previous section, our new estimates can be used to also
derive almost sure regularity results.

\begin{prop}
\label{propt} Let $H\in (0,1)$. Assume hypothesis \textnormal{(\ref{hq})}.
Let $T>0$, $t_{0}\in (0,1]$ and $s,t\in \lbrack t_{0},T]$ with $|t-s|\leq 
\frac{t_{0}}{2}$ be fixed. Then the canonical metric $\delta _{x}\left(
t-s\right) $ for $u\left( \cdot ,x\right) $ satisfies for every $x\in S^{1}$ 
\begin{equation}
c_{t_{0},T,H}|t-s|^{\alpha \wedge (2H)}\leq \delta _{x}^{2}\left( t-s\right)
\leq C_{t_{0},T,H}|t-s|^{\alpha \wedge (2H)},  \label{delta1}
\end{equation}
where $c_{t_{0},T,H}$ and $C_{t_{0},T,H}$ are positive constant depending
only on $t_{0},T$ and $H$. In particular, $u(\cdot ,x)$ is $\beta $-H\"{o}
lder continuous for any $\beta \in (0,\frac{\alpha }{2}\wedge H)$.

In particular, $u(\cdot ,x)$ is $\beta $-H\"{o}lder continuous for any $
\beta \in (0,\frac{\alpha }{2}\wedge H)$. More specifically, up to a
non-random constant, the function $r\mapsto r^{\frac{\alpha }{2}\wedge
H}\log ^{1/2}\left( 1/r\right) $ is an almost sure uniform modulus of
continuity for $u\left( \cdot ,x\right) $.
\end{prop}

\begin{proof}
The statement on almost-sure continuity is established using the arguments
described in the introduction, or simply by applying Dudley's entropy
upper bound theorem (see \cite[Theorem 2.7.1]{Khoshnevisan:02}). We detail only the proof of (\ref{delta1}), separating
the cases $H>1/2$ and $H<1/2$.

Fix $T>0$, $t_{0}\in (0,1]$ and $s,t\in \lbrack t_{0},T]$ such that $
|t-s|\leq \frac{t_{0}}{2}$. We assume without loss of generality that $s\leq
t$. Following \cite[Section 2.1]{Sarol:06}, it yields that 
\begin{equation}
\begin{split}
\delta _{x}^{2}(s,t)& =q_{0}|t-s|^{2H} \\
& +\sum_{n=1}^{+\infty }q_{n}\,{\mathrm{E}}\biggl[\biggl\{
\int_{0}^{s}(e^{-n^{2}(t-r)}-e^{-n^{2}(s-r)})\beta
_{n}^{H}(dr)+\int_{s}^{t}e^{-n^{2}(t-r)}\beta _{n}^{H}(dr)\biggr\}^{2}\biggr],
\end{split}
\label{deltatime}
\end{equation}
where $\{(\beta _{n}^{H}(t),t\geq 0)\}_{n\geq 1}$ is a sequence of fractional Brownian motions.

In order to bound the last expectation we consider two different cases:

\vskip 12pt

\noindent \textit{Case 1}: $H \geq \frac12$. In \cite[(15)]{Sarol:06} it is
proved that $\delta_x^2(s,t)$ is bounded above and below by 
\begin{equation*}  \label{eqqn}
q_0 \vert t-s \vert^{2H} + \sum_{n^2(t-s) > 1} \frac{c_H q_n}{n^{4H}}+
\sum_{n^2(t-s) \leq 1} C_H q_n \vert t-s \vert^{2H}.
\end{equation*}
Taking $q_n$ and $\alpha \in (0,1]$ from hypothesis \textnormal{(\ref{hq})},
we obtain that $\delta_x^2(s,t)$ is bounded above and below by 
\begin{equation*}
c_H (\vert t-s \vert^{2H} + \vert t-s \vert^{\alpha}).
\end{equation*}
Therefore, the upper and the lower bounds of (\ref{delta1}) follow for $H \geq \frac12$. 

\vskip 12pt

\noindent \textit{Case 2}: $H <\frac12$. We prove the upper and lower bound
of (\ref{delta1}) separately.

\vskip12pt

\noindent \textit{The upper bound}. In order to prove the upper bound of 
(\ref{delta1}), we start estimating the expectation in (\ref{deltatime}).
Using the results in Section 2.2, we have that 
\begin{equation}
{\mathrm{E}}\,\biggl[\biggl(\int_{0}^{s}\biggl(
e^{-n^{2}(t-r)}-e^{-n^{2}(s-r)}\biggr)\beta
_{n}^{H}(dr)+\int_{s}^{t}e^{-n^{2}(t-r)}\beta _{n}^{H}(dr)\biggr)^{2}\biggr]
\leq 2I_{1}+I_{2}+2I_{3},
\end{equation}
where 
\begin{equation}
\begin{split}
I_{1}& :=\int_{0}^{s}(K_{s}^{\ast
}f(r))^{2}dr,\;\;f(r)=e^{-n^{2}(t-r)}-e^{-n^{2}(s-r)}, \\
I_{2}& :=\int_{s}^{t}(K_{t}^{\ast }g(r))^{2}dr,\;\;g(r)=e^{-n^{2}(t-r)}, \\
I_{3}& :=\int_{0}^{s}(K_{t}^{\ast }g(r)-K_{s}^{\ast }g(r))^{2}dr.
\end{split}
\label{is}
\end{equation}

We start estimating $I_1$. We write 
\begin{equation}  \label{isuns}
\begin{split}
I_1 &\leq 2 \int_0^s (K(s,r) f(r))^2 dr + 2 \int_0^s \biggl( \int_r^s
(f(u)-f(r)) \frac{\partial K}{\partial u }(u,r) du \biggr)^2 dr \\
&:=2 I_{1,1}+2 I_{1,2}.
\end{split}
\end{equation}
Using Lemma \ref{lema1} and the change of variables $2n^2(s-r)=v$, we have 
\begin{equation*}
\begin{split}
I_{1,1}& \leq c_H \int_0^s (s-r)^{2H-1} r^{2H-1} ( e^{-n^2(t-r)}-
e^{-n^2(s-r)} )^2 \, dr \\
&=\frac{c_H}{n^{4H}} ( 1-e^{-n^2(t-s)})^2 \int_0^{2n^2s} \biggl(s-\frac{v}{
2n^2} \biggr)^{2H-1} v^{2H-1} e^{-v} \, dv.
\end{split}
\end{equation*}
By Lemma \ref{a1}, it yields 
\begin{equation*}
I_{1,1} \leq \frac{c_H}{n^{4H}} ( 1-e^{-n^2(t-s)})^2.
\end{equation*}
We now treat $I_{1,2}$. Using Lemma \ref{lema1} and the change of variables $
s-r=v$, $s-u=v^{\prime}$, we have 
\begin{equation*}
\begin{split}
I_{1,2}\leq c_H ( 1-e^{-n^2(t-s)})^2 \int_0^s dv \biggl(\int_0^v
dv^{\prime}(v-v^{\prime})^{H-\frac32} (e^{-n^2v^{\prime}}-e^{-n^2v})\biggr)
^2.
\end{split}
\end{equation*}
By the change of variables $v-v^{\prime}=u$, we find 
\begin{equation*}
I_{1,2}\leq c_H ( 1-e^{-n^2(t-s)})^2 \int_0^s dv \, e^{-2n^2v} \biggl( 
\int_0^{v} du \, u^{H-\frac32} (e^{n^2u}-1) \biggr)^2.
\end{equation*}
Then using \cite[Lemma 2]{Tindel:03} with $a=n^2$ and $A=H-\frac12$, we
conclude that 
\begin{equation*}
I_{1,2}\leq \frac{c_H}{n^{4H}} ( 1-e^{-n^2(t-s)})^2.
\end{equation*}
Writing $I_{1,1}$ and $I_{1,2}$ together, we get 
\begin{equation*}
I_1\leq \frac{c_H}{n^{4H}} ( 1-e^{-n^2(t-s)})^2.
\end{equation*}

We now separate the sum in (\ref{deltatime}) in two terms, as $n^2(t-s) > 1$
(tail) and $n^2(t-s) \leq 1$ (head), and take $q_n$ and $\alpha \in (0,1]$
from hypothesis \textnormal{(\ref{hq})}. Then we obtain for the tail of the
series 
\begin{equation*}
\sum_{n^2(t-s) > 1} q_n I_1 \leq c_H \sum_{n^2(t-s) > 1} n^{-2 \alpha-1}
\leq c_H \vert t-s \vert^{\alpha}.
\end{equation*}
For the head of the series, use the inequality $1-e^{-x} \leq x$, valid for
all $x\geq 0$, to get 
\begin{equation*}
\begin{split}
\sum_{n^2(t-s) \leq 1} q_n I_1 &\leq \sum_{n^2(t-s) \leq 1} q_n \frac{
c(t_0,H)}{n^{4H}} ( 1-e^{-n^2(t-s)})^{2H} ( 1-e^{-n^2(t-s)})^{2-2H} \\
&\leq c_H \vert t- s\vert^{2H} \sum_{n^2(t-s) \leq 1} n^{4H-2\alpha-1} \\
& \leq c_H \vert t- s\vert^{\alpha \wedge (2H)}.
\end{split}
\end{equation*}

We now bound $I_2$. 
\begin{equation*}
\begin{split}
I_2 &\leq 2 \int_s^t (K(t,r) g(r))^2 dr + 2 \int_s^t dr \biggl( \int_r^t du
\, (g(u)-g(r)) \frac{\partial K}{\partial u }(u,r) \biggr)^2 \\
&:=2 I_{2,1}+2 I_{2,2}.
\end{split}
\end{equation*}
Using Lemma \ref{lema1} and the change of variables $2n^2(t-r)=u$, we have 
\begin{equation*}
\begin{split}
I_{2,1} &\leq c_H \int_s^t dr \, (t-r)^{2H-1} r^{2H-1} e^{-2n^2(t-r)} \\
&=\frac{c_H}{n^{4H}} \int_0^{2n^2(t-s)} du \, \biggl(t-\frac{u}{2n^2} \biggr)
^{2H-1} u^{2H-1} e^{-u}.
\end{split}
\end{equation*}
Using Lemma \ref{a1}, we obtain for the tail of the series 
\begin{equation*}
\sum_{n^2(t-s) > 1} q_n I_{2,1} \leq c_H \sum_{n^2(t-s) > 1} n^{-2 \alpha-1}
\leq c_H \vert t-s \vert^{\alpha}.
\end{equation*}
For the head of the series, as $\vert t-s \vert \leq \frac{t_0}{2}$, we have 
\begin{equation*}
\begin{split}
\sum_{n^2(t-s) \leq 1} q_n I_{2,1} &\leq \sum_{n^2(t-s) \leq 1} q_n \frac{
c(t_0,H)}{n^{4H}} \biggl(\frac{t}{2} \biggr)^{2H-1} \int_0^{2n^2(t-s)} du \,
u^{2H-1} \\
&\leq c_H \vert t- s\vert^{2H} \sum_{n^2(t-s) \leq 1} n^{4H-2\alpha-1}.
\end{split}
\end{equation*}
This proves that $\sum_{n^{2}\left( t-s\right) \leq 1}q_{n}I_{2,1}$ is of
the same order as $\sum_{n^{2}\left( t-s\right) \leq 1}q_{n}I_{1}$ which we
calculated above to be of order $\vert t- s\vert^{\alpha \wedge (2H)}$.

We now bound $I_{2,2}$. Using Lemma \ref{lema1} and the change of variables $
t-r=v$, $t-u=v^{\prime}$, we have 
\begin{equation*}
\begin{split}
I_{2,2} \leq c_H \int_0^{t-s} dv \biggl(\int_0^v
dv^{\prime}(v-v^{\prime})^{H-\frac32} (e^{-n^2v^{\prime}}-e^{-n^2v})\biggr)
^2.
\end{split}
\end{equation*}
Using the change of variables $n^2(v-v^{\prime})=y$ and $2n^2v=x$, we find 
\begin{equation*}
I_{2,2}\leq \frac{c_H}{n^{4H}} \int_{0}^{2n^2(t-s)} dx \, e^{-x} \biggl( 
\int_0^{x/2} dy \, y^{H-\frac32} (e^{y}-1) \biggr)^2
\end{equation*}
Appealing to \cite[Lemma 2]{Tindel:03} with $a=n^2$ and $A=H-\frac12$, we
obtain for the tail of the series 
\begin{equation*}
\sum_{n^2(t-s) > 1} q_n I_{2,2} \leq c_H \sum_{n^2(t-s) > 1} n^{-2 \alpha-1}
\leq c_H \vert t-s \vert^{\alpha}.
\end{equation*}
For the head of the series, we have 
\begin{equation*}
\begin{split}
\sum_{n^2(t-s) \leq 1} q_n I_{2,2} &\leq \sum_{n^2(t-s) \leq 1} q_n \frac{c_H
}{n^{4H}} \int_{0}^{2n^2(t-s)} dx \, \biggl( \int_0^{1/2} dy y^{H-\frac32}
(e^{y}-1) \biggr)^2 \\
&\leq c_H \vert t- s\vert \sum_{n^2(t-s) \leq 1} n^{-2\alpha+1} \leq c_H
\vert t- s\vert^{\alpha}.
\end{split}
\end{equation*}

We now estimate $I_3$. 
\begin{equation*}
\begin{split}
I_3 &\leq 2 \int_0^s (K(t,r)-K(s,r))^2(g(r))^2 dr + 2 \int_0^s dr \biggl( 
\int_s^t du \, (g(u)-g(r)) \frac{\partial K}{\partial u }(u,r) \biggr)^2 \\
&:=2I_{3,1}+2I_{3,2}.
\end{split}
\end{equation*}
By Lemma \ref{lema1}, we get, for every $r<s<t$, 
\begin{align}
K\left( t,r\right) -K\left( s,r\right) & =\left( t-s\right) \int
_{0}^{1}\left\vert \frac{\partial K}{\partial u}\left( s+v\left( t-s\right)
,r\right) \right\vert dv  \notag \\
& \leq c_{H}\left( t-s\right) \int_{0}^{1}\left\vert s+v\left( t-s\right)
-r\right\vert ^{H-3/2}dv  \label{Kinc2} \\
& \leq c_{H}\left\vert \left( s-r\right) ^{H-1/2}-\left( t-r\right)
^{H-1/2}\right\vert .  \label{Kinc}
\end{align}

We now separate the evaluation of the integral in $I_{3,1}$ depending upon
whether $r$ is bigger or smaller than $s-\left( t-s\right) /2$. In the first
case, we evaluate 
\begin{equation*}
I_{3,1,1}:=\int_{s-\left( t-s\right) /2}^{s}\left( K\left( t,r\right)
-K\left( s,r\right) \right) ^{2}e^{-n^{2}\left( t-r\right) }dr.
\end{equation*}
Here, we have $s-r<\left( t-s\right) /2$ and $t-r>t-s$; therefore, using (
\ref{Kinc}), we have 
\begin{align*}
I_{3,1,1} & \leq c_{H}\int_{s-\left( t-s\right) /2}^{s}\left( \left(
1+2^{H-1/2}\right) \left( s-r\right) ^{H-1/2}\right) ^{2}e^{-n^{2}\left(
t-r\right) }dr \\
& \leq c_{H}e^{-n^{2}\left( t-s\right) }\int_{s-\left( t-s\right) /2}
^{s}\left( s-r\right) ^{2H-1}dr \\
& =c_{H}e^{-n^{2}\left( t-s\right) }\left( t-s\right) ^{2H}.
\end{align*}
For the head of the series, we find 
\begin{equation*}
\sum_{n^{2}\left( t-s\right) \leq 1}q_{n}I_{3,1,1}\leq c_{H}\left(
t-s\right) ^{2H}\sum_{n^{2}\left( t-s\right) \leq 1}n^{4H-2\alpha-1},
\end{equation*}
which is bounded above by $c_H \vert t- s\vert^{\alpha \wedge (2H)}$ while for
the tail of the series we have 
\begin{align*}
\sum_{n^{2}\left( t-s\right) > 1}q_{n}I_{3,1,1} & \leq c_{H}\left(
t-s\right) ^{2H}\sum_{n^{2}\left( t-s\right) > 1}n^{4H-2\alpha-1}
e^{-n^{2}\left( t-s\right) } \\
& \leq c_{H}\left( t-s\right) ^{2H}\int_{\left( t-s\right) ^{-1/2}
}^{\infty}e^{-x^{2}\left( t-s\right) }x^{4H-2\alpha-1}dx \\
& =c_{H}\left( t-s\right) ^{\alpha}\int_{1}^{\infty}e^{-y^{2}}
y^{4H-2\alpha-1}dy=c_{H}\left( t-s\right)^{\alpha}.
\end{align*}

Second we evaluate 
\begin{equation*}
I_{3,1,2}:=\int_{0}^{s-\left( t-s\right) /2}\left( K\left( t,r\right)
-K\left( s,r\right) \right) ^{2}e^{-n^{2}\left( t-r\right) }dr.
\end{equation*}
Here, we have $s-r>\left( t-s\right) /2$; we simply use (\ref{Kinc2}) where
an upper bound is obtained by replacing $\left\vert s+v\left( t-s\right)
-r\right\vert ^{H-3/2}$ by $\left\vert s-r\right\vert ^{H-3/2}$; the latter
can now be bounded above by $2^{3/2-H}\left\vert t-s\right\vert ^{H-3/2}$.
Thus 
\begin{align*}
I_{3,1,2} & \leq c_{H}\left\vert t-s\right\vert ^{2+2H-3}\int_{0}^{s-\left(
t-s\right) /2}e^{-n^{2}\left( t-r\right) }dr. \\
& \leq c_{H}\left\vert t-s\right\vert ^{2H-1}n^{-2}e^{-n^{2}\left(
t-s\right) }.
\end{align*}
This estimate will not help us in the case $n^{2}\left( t-s\right) \leq 1$.
In the other case, we have 
\begin{align*}
\sum_{n^{2}\left( t-s\right) > 1}q_{n}I_{3,1,2} & \leq c_{H}\left\vert
t-s\right\vert ^{2H-1}\sum_{n^{2}\left( t-s\right) > 1}n^{4H-2\alpha
-3}e^{-n^{2}\left( t-s\right) } \\
& \leq c_{H}\left\vert t-s\right\vert ^{2H-1} \int_{\left( t-s\right)
^{-1/2}}^{\infty} x^{4H-2\alpha-3}e^{-x^{2}\left( t-s\right) }dx \\
& =c_{H}\left\vert t-s\right\vert ^{2H-1}\left( t-s\right) ^{-2H+\alpha
+1}\int_{1}^{\infty}y^{4H-2\alpha-3}e^{-y^{2}}dy=c_{H}\left( t-s\right)
^{\alpha}.
\end{align*}

The third and last step of the estimation of $I_{3,1}$ is the sum for $
n^{2}\left( t-s\right) <1$ of $I_{3,1,2}$. In this case, we use (\ref{Kinc})
and obtain an upper bound by bounding $\left( s-r\right) ^{H-1/2}-\left(
t-r\right) ^{H-1/2}$ above by $c_{H}\left( t-s\right) \left( s-r\right)
^{H-3/2}$. Thus 
\begin{equation*}
I_{3,1,2}\leq c_{H}\left( t-s\right) ^{2}\int_{0}^{s-\left( t-s\right)
/2}\left( s-r\right) ^{2H-3}dr\leq c_{H}\left( t-s\right) ^{2H}.
\end{equation*}
This proves that $\sum_{n^{2}\left( t-s\right) \leq 1}q_{n}I_{3,1,2}$ is of
the same order as $\sum_{n^{2}\left( t-s\right) \leq 1}q_{n}I_{3,1,1}$ which
we calculated above to be of order $\vert t- s\vert^{\alpha \wedge (2H)}$.

We now bound $I_{3,2}$. Using Lemma \ref{lema1} and the change of variables $
s-r=v$, $s-u=v^{\prime}$, we have 
\begin{equation*}
\begin{split}
I_{3,2}&\leq c_H e^{-2n^2(t-s)} \int_0^s dv \biggl(\int_{s-t}^0
dv^{\prime}(v-v^{\prime})^{H-\frac32} (e^{-n^2v^{\prime}}-e^{-n^2v})\biggr)
^2.
\end{split}
\end{equation*}
Using the change of variables $v-v^{\prime}=u$, we find 
\begin{equation*}
I_{3,2} \leq c_H e^{-2n^2(t-s)} \int_{0}^s dv \, e^{-2n^2v} \biggl( 
\int_{v}^{v+(t-s)} du \, u^{H-\frac32} (e^{n^2u}-1) \biggr)^2.
\end{equation*}

Appealing to \cite[Lemma 2]{Tindel:03} with $a=n^2$ and $A=H-\frac12$, we
obtain for the tail of the series 
\begin{equation*}
\sum_{n^2(t-s) > 1} q_n I_{3,2} \leq c_H \sum_{n^2(t-s) > 1} n^{-2 \alpha-1}
\leq c_H \vert t-s \vert^{\alpha}.
\end{equation*}
In order to evaluate the head of the series, we separate the evaluation of
the integral in $I_{3,2}$ depending upon whether $v$ is bigger or smaller
than $t-s$, that is, 
\begin{equation*}
\begin{split}
I_{3,2} &\leq c_H \int_{0}^{s} dv \, \biggl( \int_{v}^{v+(t-s)} du \,
u^{H-\frac32} \biggr)^2 \\
&=c_H \biggl\{ \int_{0}^{t-s} dv \, \biggl( \int_{v}^{v+(t-s)} du \,
u^{H-\frac32} \biggr)^2+ \int_{t-s}^{s} dv \, \biggl( \int_{v}^{v+(t-s)} du
\, u^{H-\frac32} \biggr)^2\biggr\} \\
&\leq c_H \biggl\{ \int_{0}^{t-s} dv \, v^{2H-1}+ \int_{t-s}^s dv \,
v^{2H-3} (t-s)^2 \biggr\} \\
& \leq c_H (t-s)^{2H}.
\end{split}
\end{equation*}
Therefore, $\sum_{n^{2}\left( t-s\right) \leq 1}q_{n}I_{3,2}$ is of the same
order as $\sum_{n^{2}\left( t-s\right) \leq 1}q_{n}I_{3,1}$ which is of
order $\vert t- s\vert^{\alpha \wedge (2H)}$.

Use all the estimates above, together with (\ref{deltatime}),
to conclude that 
\begin{equation*}
\delta_x^2(s,t) \leq c_H (\vert t-s \vert^{2H} + \vert t-s
\vert^{\alpha}) \leq c^{\prime}_H \vert t-s \vert^{\alpha \wedge (2H)}.
\end{equation*}
This proves the upper bound of (\ref{delta1}) when $H<1/2$.

\vskip12pt

\noindent \textit{The lower bound}: We now estimate the lower bound of the
expectation in the case $H<1/2$. We write 
\begin{equation}
{\mathrm{E}}\,\biggl[\biggl(\int_{0}^{s}\biggl(
e^{-n^{2}(t-r)}-e^{-n^{2}(s-r)}\biggr)\beta
_{n}^{H}(dr)+\int_{s}^{t}e^{-n^{2}(t-r)}\beta _{n}^{H}(dr)\biggr)^{2}\biggr]
=I_{1}+I_{2}+I_{3}+I_{4},  \label{mino}
\end{equation}
where $I_{1}$, $I_{2}$ and $I_{3}$ are as in (\ref{is}), and 
\begin{equation}
I_{4}:=\int_{0}^{s}(K_{s}^{\ast }f(r))(K_{t}^{\ast }g(r)-K_{s}^{\ast
}g(r))dr.  \label{I4}
\end{equation}
First note that $I_{1},I_{2},I_{3}\geq 0$.

We start by finding a lower bound for $I_{1}$. We have $
I_{1}:=I_{1,1}+I_{1,2}+I_{1,3}$, where $I_{1,1}$ and $I_{1,2}$ are as in (\ref{isuns}), and 
\begin{equation*}
I_{1,3}=2\int_{0}^{s}dr\,K(s,r)f(r)\int_{r}^{s}du\,(f(u)-f(r))\frac{\partial
K}{\partial u}(u,r).
\end{equation*}
The change of variables $s-r=v,s-u=w,v-w=u^{\prime }$ gives 
\begin{equation*}
I_{1}=(1-e^{-n^{2}(t-s)})^{2}\int_{0}^{s}dv\,e^{-2n^{2}v}\biggl(
K(s,s-v)+\int_{0}^{v}du^{\prime }\,\frac{\partial K}{\partial u^{\prime }}
(u^{\prime },0)(e^{n^{2}u^{\prime }}-1)\biggr)^{2}.
\end{equation*}
Appealing to Lemma \ref{lema1} in the appendix, and the change of variables $
n^{2}u^{\prime }=u,n^{2}v=x$, we obtain
\begin{eqnarray}
I_{1} &\geq &\frac{c_{H}}{n^{4H}}(1-e^{-n^{2}(t-s)})^{2}\int_{0}^{n^{2}s}dx
\,e^{-2x}\biggl(x^{H-\frac{1}{2}}-(\frac{1}{2}-H)\int_{0}^{x}du\,u^{H-\frac{3
}{2}}(e^{u}-1)\biggr)^{2}  \notag \\
&\geq &\frac{c_{H}}{n^{4H}}(1-e^{-n^{2}(t-s)})^{2}\int_{0}^{t_{0}}dx\,e^{-2x}
\biggl(x^{H-\frac{1}{2}}-(\frac{1}{2}-H)\int_{0}^{x}du\,u^{H-\frac{3}{2}
}(e^{u}-1)\biggr)^{2}  \notag \\
&=&\frac{c_{H}}{n^{4H}}(1-e^{-n^{2}(t-s)})^{2},  \label{i1lower}
\end{eqnarray}
as the last integral is finite and positive.

Next we evaluate $I_{4}$. We write $I_{4}=I_{4,1}+I_{4,2}+I_{4,3}+I_{4,4}$,
where 
\begin{equation}
\begin{split}
I_{4,1}& =\int_{0}^{s}drK(s,r)f(r)(K(t,r)-K(s,r))g(r), \\
I_{4,2}& =\int_{0}^{s}drK(s,r)f(r)\int_{s}^{t}du(g(u)-g(r))\frac{\partial K}{
\partial u}(u,r), \\
I_{4,3}& =\int_{0}^{s}dr\int_{r}^{s}du(f(u)-f(r))\frac{\partial K}{\partial u
}(u,r)\int_{s}^{t}dv(g(v)-g(r))\frac{\partial K}{\partial v}(v,r), \\
I_{4,4}& =\int_{0}^{s}dr(K(t,r)-K(s,r))g(r)\int_{r}^{s}du(f(u)-f(r))\frac{
\partial K}{\partial u}(u,r).
\end{split}
\label{I4s}
\end{equation}
Now, note that $I_{4,1},I_{4,2}\geq 0$ but $I_{4,3},I_{4,4}\leq 0$.

We claim that, for some subset $S_{K}\subset \mathbb{N}$, 
\begin{align}
& \sum_{n\in S_{K}}q_{n}I_{1}>2\sum_{n\in S_{K}}q_{n}\left\vert
I_{4,3}\right\vert ,  \label{I1ge2I43} \\
& \sum_{n\in S_{K}}q_{n}I_{1}>4\sum_{n\in S_{K}}q_{n}\left\vert
I_{4,4}\right\vert ,  \label{I1ge2I44}
\end{align}
where $q_{n}$ and $\alpha \in (0,1]$ are as in hypothesis \textnormal{(\ref
{hq})} and $S_{K}:=\left\{ n\in \mathbb{N}:n^{2}\left( t-s\right) >K\right\} 
$ for some (large) constant $K\geq 1$ which will be chosen later.

Assume (\ref{I1ge2I43}) and (\ref{I1ge2I44}) proved. We write, using (\ref
{i1lower}),
\begin{equation}
\sum_{n\in S_{K}}q_{n}I_{1}\geq c_{H}(1-e^{-1})^{2}\int_{2\sqrt{K}/\sqrt{t-s}
}^{\infty }dx\,x^{-2\alpha -1}:=c_{\alpha ,H}^{1}(t-s)^{\alpha }{K}^{-\alpha
}.  \label{lowi1}
\end{equation}

Because $I_2, I_3, I_{4,1}, I_{4,2} \geq 0$ and using (\ref{I1ge2I43}), (\ref
{I1ge2I44}) and (\ref{lowi1}), we find 
\begin{align*}
\sum_{n\in \mathbb{N}} q_n (I_1+I_2+I_3+I_4) &\geq \sum_{n\in S_K} q_n I_{1}
-\sum _{n\in S_K} q_n \vert I_{4,3} \vert- \sum _{n\in S_K} q_n \vert I_{4,4} \vert \\
&\geq\frac{1}{4} \sum_{n\in S_K} q_n I_{1} \\
& \geq c_{\alpha, H, K} \left( t-s\right) ^{\alpha}.
\end{align*}

Therefore, by (\ref{deltatime}) and (\ref{mino}), we conclude that 
\begin{equation*}
\delta_x^2(s,t) \geq q_0 \vert t-s \vert^{2H} + c_H \vert t-s \vert^{\alpha} \geq c_H^{\prime} \vert t-s \vert^{\alpha \wedge (2H)}.
\end{equation*}
This proves the lower
bound of (\ref{delta1}) when $H < \frac12$.

\vskip 12pt

We finally prove (\ref{I1ge2I43}) and (\ref{I1ge2I44}).
\begin{proof}[Proof of \textnormal{(\protect\ref{I1ge2I43})}]
Using Lemma \ref{lema1} and the change of variables $s-r=r^{\prime}$, $
s-u=u^{\prime}$, $s-v=v^{\prime}$, $r^{\prime}-u^{\prime}=u^{\prime\prime}$, 
$r^{\prime}-v^{\prime}=v^{\prime\prime}$, $n^2u^{\prime\prime}=x$, $
n^2v^{\prime\prime}=v$, we find 
\begin{equation}
\begin{split}
\left\vert I_{4,3}\right\vert &\leq c_H \left( 1-e^{-n^2 (t-s)}\right)
e^{-n^2 (t-s)} n^{-4H} \int_{0}^{n^{2}s} dxe^{-2x} \left( \int_{0}^{x} du \,
u^{H-3/2}\left( e^{u}-1\right) \right) \\
&\qquad \qquad \qquad \qquad \qquad \qquad \times \left( \int_{x}^{x+n^2
(t-s) } dv \, v^{H-3/2}\left( e^{v}-1\right) \right).  \label{I43UB0}
\end{split}
\end{equation}

Note that with the exception of the factor $e^{-n^2(t-s)}$ in $\left\vert
I_{4,3} \right\vert $, the combination of all the terms in $I_{1}$ and $
I_{4,3}$ are in fact largely similar, which makes this portion of the proof
quite delicate, and in particular, to exploit the factor $e^{-n^2 (t-s)}$,
we must restrict the values of $n^2(t-s)$ to being relatively large, which
explains the choice of $S_K$ above.

Our strategy is to bound the sum over $n\in S_K$ of $q_n \left\vert I_{4,3}
\right\vert $ above as tightly as possible by performing a \textquotedblleft
Fubini\textquotedblright, dragging the sum over $n$ all the way inside the
expression for $\sum_{S_{K}} q_n \left\vert I_{4,3}\right\vert $, and
evaluating it first using some Gaussian estimates. That these Gaussian
estimates work has to do with the precise eigenvalue structure of the
Laplacian, not with the Gaussian property of the driving noise.

We proved in (\ref{lowi1}) that the contribution of $I_{1}$ is bounded below
by an expression of the form $c^1_{\alpha, H} \left( t-s\right)
^{\alpha}K^{-\alpha}$, where $c^1_{\alpha, H}$ depends only on $\alpha$ and $
H$. We will now show that 
\begin{equation}  \label{new}
\sum_{n\in S_K} q_n \vert I_{4,3} \vert \leq c_{\alpha, H}^2 \left(
t-s\right) ^{\alpha} K^{-\beta}.
\end{equation}
for some $\beta>\alpha$, where $c_{\alpha, H}^{2}$ depends again only on $H$
and $\alpha$. Even if $c_{\alpha, H}^{2}$ is much larger than $
c_{\alpha,H}^1 $, one only needs to choose $K\geq(2 c_{\alpha,H}^{2}/c^1_{H,
\alpha})^{1/\left( \beta-\alpha\right) } $ to guarantee that the
contribution of $I_1$ exceeds twice the absolute value of the contribution
of $I_{4,3}$ as announced in (\ref{I1ge2I43}), which implies that even
though the latter is negative, the sum of the two exceeds $(c^1_{\alpha,
H}/2)\left( t-s\right) ^{\alpha} K^{-\alpha}$, i.e. for some $K$ depending
only on $H$ and $\alpha$.

First, for fixed $x$, we perform the announced Fubini, which means that,
instead of having the integration and summation limits for $n$ and $v$ as $n>
\sqrt{K/\left( t-s\right) }$ first and $x<v<x+n^{2}\left( t-s\right) $ next,
we get instead $x<v<\infty $ and 
\begin{align*}
n& >\max \left\{ \sqrt{K/\left( t-s\right) },\sqrt{\left( v-x\right) /t-s}
\right\} \\
& =\left( t-s\right) ^{-1/2}\sqrt{\left( v-x\right) \vee K}.
\end{align*}
Therefore, bounding $(1-e^{-n^{2}(t-s)})$ by $1$, and $n^{2}s$ by $\infty $,
we have 
\begin{equation}
\begin{split}
\sum_{n\in S_{K}}q_{n}\left\vert I_{4,3}\right\vert & \leq
c_{H}\int_{0}^{\infty }dx\,e^{-2x}\left( \int_{0}^{x}du\,u^{H-3/2}\left(
e^{u}-1\right) \right) \\
& \qquad \times \left( \int_{x}^{\infty }dv\,v^{H-3/2}\left( e^{v}-1\right)
S\left( K,v-x,t-s\right) \right) ,
\end{split}
\label{I43UB1}
\end{equation}
where the term $S\left( K,v-x,t-s\right) $ is defined by a series which we
compare to a Gaussian integral as follows 
\begin{align*}
S\left( K,v-x,t-s\right) & :=\sum_{n>(t-s)^{-1/2}\sqrt{(v-x)\vee K}
}n^{-2\alpha -1}e^{-n^{2}(t-s)} \\
& \leq \int_{y\geq (t-s)^{-1/2}\sqrt{(v-x)\vee K}}^{\infty }dy\,y^{-2\alpha
-1}e^{-y^{2}(t-s)}.
\end{align*}

Using the change of variable $w^{2}=\left( t-s\right) y^{2}$, we have 
\begin{align*}
S\left( K,v-x,t-s\right) & \leq (t-s)^{\alpha} \int_{ \sqrt{ (v-x)\vee K}
}^{\infty} dw \, w^{-2\alpha-1}e^{-w^{2}} \\
& \leq (t-s)^{\alpha} \left((v-x)\vee K \right) ^{-\alpha-1/2} \int_{\sqrt{
(v-x)\vee K}}^{\infty} dw \, e^{-w^{2}}.
\end{align*}
Now, using the classical Gaussian tail estimate $\int_{A}^{\infty} dw \,
e^{-w^{2} }\leq 2^{-1} A^{-1}e^{-A^{2}},$ we get 
\begin{equation}
S\left( K,v-x,t-s\right) \leq 2^{-1} (t-s)^{\alpha}\left( (v-x)\vee K\right)
^{-\alpha-1}e^{-(v-x) \vee K}.  \label{SKvt}
\end{equation}

Combining (\ref{I43UB1}) and (\ref{SKvt}) we have immediately 
\begin{align}
& \sum_{n\in S_{K}}q_{n}\left\vert I_{4,3}\right\vert \leq c_{H}\left(
t-s\right) ^{\alpha }\int_{0}^{\infty }dx\,e^{-2x}\left(
\int_{0}^{x}u^{H-3/2}\left( e^{u}-1\right) du\right)  \notag \\
& \qquad \qquad \times \left( \int_{x}^{\infty }dvv^{H-3/2}\left(
e^{v}-1\right) \left( \left( v-x\right) \vee K\right) ^{-\alpha
-1}e^{-\left( v-x\right) \vee K}\right)  \notag \\
& =c_{H}\left( t-s\right) ^{\alpha }e^{-K}K^{-\alpha -1}\int_{0}^{\infty
}dx\,e^{-2x}\left( \int_{0}^{x}du\,u^{H-3/2}\left( e^{u}-1\right) \right)
\label{I43detail3} \\
& \qquad \qquad \times \left( \int_{x}^{x+K}dv\,v^{H-3/2}\left(
e^{v}-1\right) \right)  \notag \\
& \qquad +c_{H}\left( t-s\right) ^{\alpha }\int_{0}^{\infty
}dx\,e^{-2x}\left( \int_{0}^{x}u^{H-3/2}du\,\left( e^{u}-1\right) \right)
\label{I43detail4} \\
& \qquad \qquad \times \left( \int_{x+K}^{\infty }dv\,v^{H-3/2}\left(
e^{v}-1\right) \left( v-x\right) ^{-\alpha -1}e^{-\left( v-x\right) }\right)
.  \notag
\end{align}

We separate the last expression into various terms. We will calculate first
the term in line (\ref{I43detail3}) by separating the $x$-integration over $
x\in\lbrack0,K]$ and $x\in(K,\infty)$, which we denote by $J_{4,3,1}$ and $
J_{4,3,2}$, respectively. The term in line (\ref{I43detail4}), which we
denote by $J_{4,3,2}$, can be dealt with more directly. We now perform these
evaluations.

\vskip 12pt

\emph{Term }$J_{4,3,1}$. We write 
\begin{align*}
J_{4,3,1} & :=c_{H} \left( t-s\right) ^{\alpha}e^{-K}K^{-\alpha-1}
\int_{0}^{K}dx \, e^{-2x}\left( \int_{0}^{x} du \, u^{H-3/2}\left(
e^{u}-1\right) \right) \\
& \qquad \qquad \times \left( \int_{x}^{x+K} dv \, v^{H-3/2}\left(
e^{v}-1\right)\right) \\
& \leq c_{H}\left( t-s\right) ^{\alpha}e^{-K}K^{-\alpha-1}\int
_{0}^{\infty}dx \, e^{-2x}\left( \int_{0}^{x} du \, u^{H-3/2}\left(
e^{u}-1\right) \right) \\
& \qquad \qquad \times \left( c_{H}+\int_{1}^{2K}dv \, v^{H-3/2}\left(
e^{v}-1\right) \right).
\end{align*}
Now, integrating by parts, we get 
\begin{equation*}
\int_{1}^{2K} dv \, v^{H-3/2}\left( e^{v}-1\right) \leq c_{H} e^{K}
K^{H+1/2}.
\end{equation*}
The last two estimates imply immediately that 
\begin{equation*}
J_{4,3,1}\leq c_{H}\left( t-s\right) ^{\alpha} K^{-\alpha +H-1/2},
\end{equation*}
which proves the contribution of $J_{4,3,1}$ in (\ref{new}).

\vskip 12pt

\emph{Term }$J_{4,3,2}$. We write
\begin{align*}
J_{4,3,2} & :=c_{H}\left( t-s\right) ^{\alpha}e^{-K}K^{-\alpha-1}
\int_{K}^{\infty}dx \, e^{-2x}\left( \int_{0}^{x} du \, u^{H-3/2}\left(
e^{u}-1\right) \right) \\
& \qquad \qquad \times \left( \int_{x}^{x+K} dv \, v^{H-3/2}\left(
e^{v}-1\right) \right) \\
& \leq c_{H}\left( t-s\right) ^{\alpha}e^{-K}K^{-\alpha-1}\int
_{K}^{\infty}dx \, e^{-2x}\left( \int_{0}^{x} du \, u^{H-3/2}\left( e^{u}
-1\right) \right) \\
& \qquad \qquad \times x^{H-3/2}\left( e^{x+K}-e^{x}\right) \\
& \leq c_{H}\left( t-s\right) ^{\alpha}K^{-\alpha-1}\int_{K}^{\infty }dx \,
e^{-x} x^{H-3/2} \left( \int_{0}^{x} du \, u^{H-3/2}\left( e^{u}-1\right)\right) \\
& \leq c_{H}\left( t-s\right) ^{\alpha}K^{-\alpha-1}\int_{K} ^{\infty}dx \,
e^{-x} x^{H-3/2} \left( c_H+e^{x}\int_{1}^{x} du \, u^{H-3/2}\right) \\
& \leq c_{H}\left( t-s\right) ^{\alpha}K^{-\alpha-1}\left(
K^{H-3/2}e^{-K}+K^{H-1/2}\right) \\
& \leq c_{H}\left( t-s\right) ^{\alpha}K^{-\alpha+H-3/2}.
\end{align*}
which proves the contribution of $J_{4,3,2}$ in (\ref{new}).

\vskip 12pt

\emph{Term }$J_{4,3,3}$. The last part of the estimation is that of 
\begin{align*}
J_{4,3,3} & :=c_{H}\left( t-s\right) ^{\alpha}\int_{0}^{\infty }dx \,
e^{-2x}\left( \int_{0}^{x}u^{H-3/2}\left( e^{u}-1\right) du\right) \\
& \qquad \qquad \times \left( \int_{x+K}^{\infty} dv \, v^{H-3/2}\left(
e^{v}-1\right) \left( v-x\right) ^{-\alpha-1}e^{-\left( v-x\right) }\right)
\\
& \leq c_{H}\left( t-s\right) ^{\alpha}K^{H-3/2}\int_{0}^{\infty }dx \,
e^{-x}\cdot\left( \int_{0}^{x} du \, u^{H-3/2}\left( e^{u}-1\right) \right)
\\
& \qquad \qquad \times \left( \int_{x+K}^{\infty} dv \, \left( v-x\right)
^{-\alpha-1} \right) \\
& =c_{\alpha, H}\left( t-s\right) ^{\alpha}K^{-\alpha+H-3/2}\int_{0}^{\infty
}du \, u^{H-3/2}\left( e^{u}-1\right) \left( \int_{u}^{\infty}dx \,
e^{-x}\right) \\
& =c_{\alpha}\left( t-s\right) ^{\alpha}K^{-\alpha+H-3/2}\int_{0}^{\infty
}u^{H-3/2}\left( e^{u}-1\right) e^{-u}du \\
& =c_{\alpha}\left( t-s\right) ^{\alpha}K^{-\alpha+H-3/2}\left[
c_{H}+\int_{1}^{\infty}u^{H-3/2}du\right] \\
& =c_{\alpha,H}\left( t-s\right) ^{\alpha}K^{-\alpha+H-3/2}.
\end{align*}
Therefore, (\ref{new}) holds taking $\beta=\alpha+1/2-H$ which is greater
than $\alpha$ as $H<1/2$.

The proof of (\ref{I1ge2I43}) is now finished.
\end{proof}

\begin{proof}[Proof of \textnormal{(\protect\ref{I1ge2I44})}]
By (\ref{Kinc2}) and Lemma \ref{lema1}, we have 
\begin{equation*}
\vert I_{4,4} \vert \leq c_H (t-s) \int_0^s dr \,(s-r)^{H-\frac32} g(r)
\int_r^s du \, (u-r)^{H-\frac32} (f(u)-f(r)).
\end{equation*}
Using the change of variables $s-r=r^{\prime}, s-u=u^{\prime},
r^{\prime}-u^{\prime}=v, n^2v=u, n^2r^{\prime}=x$, we get 
\begin{equation*}
\vert I_{4,4} \vert \leq \frac{c_H}{n^{4H-2}} (t-s) e^{-n^2 (t-s)}
(1-e^{-n^2 (t-s)} ) \int_0^{n^2s} dx \, x^{H-\frac32} e^{-2x} \int_0^x du \,
u^{H-\frac32} (e^u-1).
\end{equation*}
Bounding $(1-e^{-n^2 (t-s)} )$ by $1$ and $n^2s$ by $\infty$, we get 
\begin{equation*}
\vert I_{4,4} \vert \leq \frac{c_H}{n^{4H-2}} (t-s) e^{-n^2 (t-s)}.
\end{equation*}
We will now proceed as in the proof of (\ref{I1ge2I43}); that is we will
prove that there exists a constant $c^3_{H}$ depending only on $H$ such that 
\begin{equation}  \label{new2}
\sum_{n\in S_K} q_n \vert I_{4,4} \vert \leq c_{H}^3 \left( t-s\right)
^{\alpha} K^{-\beta},
\end{equation}
for some $\beta>\alpha$. It then suffices to choose $K\geq(4
c_{H}^{3}/c^1_{H, \alpha})^{1/\left( \beta-\alpha\right) } $ to get (\ref
{I1ge2I44}).

We now prove (\ref{new2}). We write 
\begin{equation*}
\begin{split}
\sum_{n\in S_K} q_n \vert I_{4,4} \vert &\leq c_H (t-s) \int_{\sqrt{K/(t-s)}
}^{\infty} dx \, x^{-2\alpha+1} e^{-x^2 (t-s)} \\
& = c_H (t-s)^{\alpha} \int_{\sqrt{K}}^{\infty} dy \, y^{-2\alpha+1} e^{-y^2}
\\
& \leq c_H (t-s)^{\alpha} K^{-\alpha} 2^{-1} \int_{\sqrt{K}}^{\infty} dy \,
2y e^{-y^2} \\
&\leq c_H (t-s)^{\alpha} K^{-(\alpha+1)},
\end{split}
\end{equation*}
which proves (\ref{new2}) taking $\beta=\alpha +1$ and concludes the proof
of (\ref{I1ge2I44}).
\end{proof}

This finishes the proof of the entire proposition.
\end{proof}

\section{Gaussian upper bound for the bivariate density}

We denote by $p_{t,x;s,y}(\cdot\,,\cdot)$ the (Gaussian) probability density
function of the random vector $(u(t\,,x)\,,u(s\,,y))$ for all $s,t>0$ and $
x,y\in S^1$ such that $(t\,,x) \neq (s\,,y)$.

For every fixed real number $0<\alpha \leq 1$ we consider the metric 
\begin{equation}
\mathbf{\Delta }((t,x);(s,y))=|x-y|^{2\alpha }+|t-s|^{\alpha \wedge (2H)}.
\label{del}
\end{equation}

In this section we establish an upper bound of Gaussian type for the
bivariate density $p_{t,x;s,y}(\cdot \,,\cdot )$ in terms of the metric (\ref
{del}). This will be one of the key results in order to show the lower bound
of Theorem \ref{t1}. The estimates obtained in the previous section to prove
space and time regularity are nearly sufficient to obtain the results in
this section. The following further improvement is needed, which deals with
precise joint regularity (see \cite[(4.11)]{Dalang:05} for the space-time white noise case).
\begin{lem}
\label{lemmatex}Assume hypothesis \textnormal{(\ref{hq})}. Fix $t_{0},T>0$.
Then there exists $c_{H}>0$ such that for any $s,t\in \lbrack t_{0},T]$, $
x,y\in S^{1}$, with $(t,x)$ is sufficiently near $(s,y)$, and $i=1,...,d$,
\begin{equation}
c_{H}^{-1}\mathbf{\Delta }((t\,,x)\,;(s\,,y))\leq {\mathrm{E}}\left[
(u_{i}(t\,,x)-u_{i}(s\,,y))^{2}\right] \leq c_{H}\mathbf{\Delta }
((t\,,x)\,;(s\,,y)).  \label{incr}
\end{equation}
\end{lem}

\begin{proof}
The upper bound in \eqref{incr} is a consequence of the upper bounds of
Corollary \ref{corx} and Proposition \ref{propt}, and the following
inequality 
\begin{equation*}
\E\left[ (u_{i}(t\,,x)-u_{i}(s\,,y))^{2}\right] \leq 2\{\E
\left[ (u_{i}(t\,,x)-u_{i}(s\,,x))^{2}\right] +\E\left[
(u_{i}(s\,,x)-u_{i}(s\,,y))^{2}\right] \}.
\end{equation*}

We now proceed to the proof of the lower bound in \eqref{incr}. By Corollary 
\ref{corx}, there exist $c_{1}, c_2>0$ such that for all $t\in \lbrack t_{0},T]$
, $x,y\in S^{1}$, with $x$ is sufficiently near $y$, and $i=1,...,d$, 
\begin{equation}
c_{1}|x-y|^{2\alpha }\leq \E\left[ (u_{i}(t\,,x)-u_{i}(t\,,y))^{2}%
\right] \leq c_{2}|x-y|^{2\alpha }.  \label{incr2a}
\end{equation}
Moreover, Proposition \ref{propt} assures the existence of $c_{3}, c_{4}>0$
such that that for any $s,t\in \lbrack t_{0},T]$, $x\in S^{1}$, with $t$ is sufficiently near $s$, and $i=1,...,d$, 
\begin{equation}
c_{3}\left\vert t-s\right\vert ^{\alpha \wedge (2H)}\leq \E\left[
(u_{i}(t\,,x)-u_{i}(t\,,y))^{2}\right] \leq c_{4}\left\vert t-s\right\vert
^{\alpha \wedge (2H)}.  \label{incr2b}
\end{equation}

Let us now consider two different cases. \vskip12pt

\noindent \textit{Case 1}: $|t-s|^{\alpha \wedge (2H)}<\frac{c_{1}}{4c_{4}}
|x-y|^{2\alpha }$. Appealing to the lower bound in \eqref{incr2a} and the
upper bound in \eqref{incr2b},
\begin{equation*}
\begin{split}
\E\left[ (u_{i}(t\,,x)-u_{i}(s\,,y))^{2}\right] & =\E\left[
(u_{i}(t\,,x)-u_{i}(t\,,y)+u_{i}(t\,,y)-u_{i}(s\,,y))^{2}\right] \\
& \geq \frac{1}{2}\E\left[ (u_{i}(t\,,x)-u_{i}(t\,,y))^{2}\right] -
\E\left[ (u_{i}(t\,,y)-u_{i}(s\,,y))^{2}\right] \\
& \geq \frac{1}{2}c_{1}|x-y|^{2\alpha }-c_{4}|t-s|^{\alpha \wedge (2H)}.
\end{split}
\end{equation*}
Because of the inequality that defines this Case 1, this is bounded below by 
\begin{equation*}
\begin{split}
\frac{c_{1}}{2}|x-y|^{2\alpha }-\frac{c_{1}}{4}|x-y|^{2\alpha }& =\frac{c_{1}
}{4}|x-y|^{2\alpha } \\
& \geq \frac{c_{1}}{8}|x-y|^{2\alpha }+\frac{c_{1}}{8}\frac{4c_{4}}{c_{1}}
|t-s|^{\alpha \wedge (2H)} \\
& \geq \min \left( \frac{c_{1}}{8},\frac{c_{4}}{2}\right) \mathbf{\Delta }
((t\,,x)\,;(s\,,y)).
\end{split}
\end{equation*}
This completes the proof of the lower bound in (\ref{incr}) in Case 1.\bigskip

\noindent \textit{Case 2}: $|t-s|^{\alpha \wedge (2H)}>\frac{4c_{2}}{c_{3}}
|x-y|^{2\alpha }$. The proof of this portion is identical to Case 1, by
using the upper bound in \eqref{incr2a} and the lower bound in \eqref{incr2b}, and writing 
\begin{eqnarray*}
\E\left[ (u_{i}(t\,,x)-u_{i}(s\,,y))^{2}\right] &=&\E\left[
(u_{i}(t\,,x)-u_{i}(t\,,y)+u_{i}(t\,,y)-u_{i}(s\,,y))^{2}\right] \\
&\geq &\frac{1}{2}\E\left[ \left( u_{i}\left( t,x\right)
-u_{i}\left( s,x\right) \right) ^{2}\right] -\E\left[ \left(
u_{i}\left( s,x\right) -u_{i}\left( s,y\right) \right) ^{2}\right]
\end{eqnarray*}
which yields the lower bound $\min \left( \frac{c_{3}}{8},\frac{c_{2}}{2}
\right) \mathbf{\Delta }((t\,,x)\,;(s\,,y))$.  This completes the proof of
Case 2.\bigskip

\noindent \textit{Case 3}: $\frac{4c_{2}}{c_{3}}|x-y|^{2\alpha }\geq
|t-s|^{\alpha \wedge (2H)}\geq \frac{c_{1}}{4c_{4}}|x-y|^{2\alpha }$. Note
that it suffices to prove that, 
\begin{equation}
\E\left[ (u_{i}(t\,,x)-u_{i}(s\,,y))^{2}\right] \geq c|t-s|^{\alpha
\wedge (2H)}.  \label{claim1}
\end{equation}
Indeed, because of the lower bound inequality that defines this Case 3, this
is bounded below by 
\begin{equation*}
\frac{c}{2}|t-s|^{\alpha \wedge (2H)}+\frac{c}{2}\frac{c_{1}}{4c_{4}}
|x-y|^{2\alpha }\geq c^{\prime }\mathbf{\Delta }((t\,,x)\,;(s\,,y)),
\end{equation*}
which proves the lower bound in (\ref{incr}) in this Case 1, provided that 
\eqref{claim1} is proved.

\vskip 12pt

\noindent \textit{Proof of \eqref{claim1}}. 
We write 
\begin{equation*}
\E\left[ (u_{i}(t\,,x)-u_{i}(s\,,y))^{2}\right] =q_{0}|t-s|^{2H}+
\sum_{n=1}^{\infty }q_{n}\{\mathcal{W}_{1}+\mathcal{W}_{2}\},
\end{equation*}
where 
\begin{equation*}
\begin{split}
\mathcal{W}_{1}& =\E\biggl[\biggl\{\int_{0}^{s}(\cos
(nx)\,e^{-n^{2}(t-r)}-\cos (ny)\,e^{-n^{2}(s-r)})\,\beta _{n}(dr) \\
& \qquad \qquad \qquad \qquad \qquad \qquad +\int_{s}^{t}\cos
(nx)\,e^{-n^{2}(t-r)}\,\beta _{n}(dr)\biggr\}^{2}\biggr], \\
\mathcal{W}_{2}& =\E\biggl[\biggl\{\int_{0}^{s}(\sin
(nx)\,e^{-n^{2}(t-r)}-\sin (ny)\,e^{-n^{2}(s-r)})
\,\beta _{n}^{\prime }(dr) \\
& \qquad \qquad \qquad \qquad \qquad \qquad +\int_{s}^{t}\sin
(nx)\,e^{-n^{2}(t-r)}\,\beta _{n}^{\prime }(dr)\biggr\}^{2}\biggr],
\end{split}
\end{equation*}
where $\{\beta _{n}\}_{n\in \mathbb{N}}$ and $\{\beta _{n}^{\prime }\}_{n\in 
\mathbb{N}}$ are independent standard fractional Brownian motions.

Now, because the further calculations use fractional stochastic calculus we need to
consider the two different cases, namely $H<\frac12$ and $H \geq \frac12$.

\vskip 12pt
\noindent {\it Case $H \geq \frac12$}. If $H<\alpha/2$, because $\E\left[ (u_{i}(t\,,x)-u_{i}(s\,,y))^{2}\right] \geq q_{0}|t-s|^{2H}$, \eqref{claim1} follows directly. Therefore, we assume that $H > \alpha/2$. In this case, note that   \eqref{claim1} is proved in \cite{Sarol:06} for the case $x=y$. 

Straightforward computations using (\ref{isometry}) give 
\begin{align*}
&\E\left[ (u_{i}(t\,,x)-u_{i}(s\,,y))^{2}\right] \\
&\qquad =q_{0}|t-s|^{2H}+
\sum_{n=1}^{\infty }q_{n} 
\biggl\{ \biggl(e^{-2n^2t}+e^{-2n^2s}-2\cos(n\vert x -y \vert) e^{-n^2(t+s)}\biggr)  I_1 \\
& \qquad \qquad \qquad +e^{-n^2 I_2} I_2 +2 e^{-n^2t} \biggl( e^{-n^2t}-\cos(n \vert x-y \vert) e^{-n^2s} \biggr) I_3 \biggr\} \\
& \qquad \geq q_{0}|t-s|^{2H}+
\sum_{n=1}^{\infty }q_{n} 
\biggl\{ (e^{-n^2t}-e^{-n^2s})^2  I_1 +e^{-n^2 I_2} I_2 +2 e^{-n^2t} ( e^{-n^2t}-e^{-n^2s}) I_3 \biggr\},
\end{align*}
where
\begin{align*}
I_1 & =\int_0^s dw \int_0^s dv \, e^{n^2 (w+v)}  \vert w-v \vert^{2H-2}, \\
I_2 & =\int_s^t dw \int_s^t dv \, e^{n^2 (w+v)}  \vert w-v \vert^{2H-2}, \\
I_3 & =\int_0^s dw \int_s^t dv \, e^{n^2 (w+v)}  \vert w-v \vert^{2H-2}.
\end{align*}
Hence, using the results of \cite[Section 2.1 and (17)]{Sarol:06} and \eqref{hq}, it follows that
\begin{align*}
\E\left[ (u_{i}(t\,,x)-u_{i}(s\,,y))^{2}\right] &\geq q_{0}|t-s|^{2H}+
c_H (t-s)^{2H} \sum_{n^2(t-s) \leq 1} q_{n} \\ 
&\geq c_H (t-s)^{\alpha}.
\end{align*}
This proves \eqref{claim1} when $H \geq \frac12$.

\vskip 12pt
\noindent {\it Case $H<\frac12$}. It is elementary to see that by \eqref{mino}, $\mathcal{W}_{1}\geq \tilde{I}
_{1}+\tilde{I}_{4}$, where $\tilde{I}_{1}$ and $\tilde{I}_{4}$ are defined,
respectively, as $I_{1}$ and $I_{4}$ in the previous section (see (\ref{is})
and (\ref{I4})), but replacing $f$ and $g$ by 
\begin{equation*}
\tilde{f}(r)=\cos (nx)\,e^{-n^{2}(t-r)}-\cos (ny)\,e^{-n^{2}(s-r)},\qquad 
\tilde{g}(r)=\cos (nx)\,e^{-n^{2}(t-r)}.
\end{equation*}
Similarly, $\mathcal{W}_{2}\geq \bar{I}_{1}+\bar{I}_{4}$, where $\bar{I}_{1}$
and $\bar{I}_{4}$ are defined, respectively, as $I_{1}$ and $I_{4}$ but
replacing $f$ and $g$ by 
\begin{equation*}
\bar{f}(r)=\sin (nx)\,e^{-n^{2}(t-r)}-\sin (ny)\,e^{-n^{2}(s-r)},\qquad \bar{
g}(r)=\sin (nx)\,e^{-n^{2}(t-r)}.
\end{equation*}
Therefore, the proof of \eqref{claim1} when $H<\frac12$ is
similar to the control of $I_{1}$ from below by $\left\vert I_{4}\right\vert 
$ in the previous section; yet it is less delicate, because the hardest
estimates we will need to use are one which were already obtained therein.
Indeed, proceeding as in \eqref{i1lower}, we find
\begin{eqnarray}
\tilde{I}_{1}+\bar{I}_{1} &\geq &\frac{c_{H}}{n^{4H}}\{(\cos
(nx)\,e^{-n^{2}(t-s)}-\cos (ny))^{2}+(\sin (nx)\,e^{-n^{2}(t-s)}-\sin
(ny))^{2}\}  \notag \\
&=&\frac{c_{H}}{n^{4H}}\{e^{-2n^{2}(t-s)}+1-2\cos (n|x-y|)\,e^{-n^{2}(t-s)}\}
\notag \\
&\geq &\frac{c_{H}}{n^{4H}}(1-e^{-n^{2}(t-s)})^{2}.  \label{IIlb}\end{eqnarray}
Here we see that the case where $x=y$ is the worst case, in the sense that
the lower bound (\ref{i1lower}) obtained for $I_{1}$ is a lower bound for
all $\tilde{I}_{1}+\bar{I}_{1}$ uniformly in $t,x,s,y$.

Moreover, simple calculations yield very similar formulas for the four terms
in $\tilde{I}_{4}+\bar{I}_{4}$ as we had found for $I_{4}$ itself in (\ref
{I4s}); namely we have
\begin{eqnarray*}
\tilde{I}_{4,1}+\bar{I}_{4,1} &=&\int_{0}^{s}drK(s,r)h(r)(K(t,r)-K(s,r))g(r),
\\
\tilde{I}_{4,2}+\bar{I}_{4,2}
&=&\int_{0}^{s}drK(s,r)h(r)\int_{s}^{t}du(g(u)-g(r))\frac{\partial K}{
\partial u}(u,r), \\
\tilde{I}_{4,3}+\bar{I}_{4,3} &=&\int_{0}^{s}dr\int_{r}^{s}du(h(u)-h(r))
\frac{\partial K}{\partial u}(u,r)\int_{s}^{t}dv(g(v)-g(r))\frac{\partial K}{
\partial v}(v,r), \\
\tilde{I}_{4,4}+\bar{I}_{4,4}
&=&\int_{0}^{s}dr(K(t,r)-K(s,r))g(r)\int_{r}^{s}du(h(u)-h(r))\frac{\partial K
}{\partial u}(u,r).
\end{eqnarray*}
where 
\begin{eqnarray}
h(r) &=&e^{-n^{2}(t-r)}-\cos (n|x-y|)\,e^{-n^{2}(s-r)}  \notag \\
&=&e^{-n^{2}\left( s-r\right) }\left( e^{-n^{2}\left( t-s\right) }-\cos
(n|x-y|)\right)  \notag \\
&=:& e^{-n^{2}\left( s-r\right) }h_{s,t,x,y}.  \label{hstxy}
\end{eqnarray}
In other words, for each $j=1,2,3,4$, the formula for $\tilde{I}_{4,j}+\bar{I
}_{4,j}$ is identical to that of $I_{4,j}$, with $f$ replaced by $h$. Also
recall that 
\begin{equation*}
f\left( r\right) =e^{-n^{2}\left( s-r\right) }\left( e^{-n^{2}\left(
t-s\right) }-1\right) =e^{-n^{2}\left( s-r\right) }h_{s,t,x,x}.
\end{equation*}
We see here that $f$ is always negative, while it is much more difficult to
control the sign of $h$. Luckily, for any $r$, the sign of $h\left( r\right) 
$ is the sign of the fixed coefficient $h_{s,t,x,y}$ defined in (\ref
{hstxy}). When $h_{s,t,x,y}$ is negative, we will be able to use
calculations from the previous section directly. When $h_{s,t,x,y}$ is
non-negative, we will instead compare $\tilde{I}_{1}+\bar{I}_{1}$ with $
\left\vert \tilde{I}_{4,1}+\bar{I}_{4,1} \right\vert $ and $\left\vert \tilde{
I}_{4,2}+\bar{I}_{4,2} \right\vert .$\bigskip

\emph{Case }$h_{s,t,x,y}<0$. Note that, in this case, $\tilde{I}_{4,1}+\bar{I}
_{4,1}>0$ and $\tilde{I}_{4,2}+\bar{I}_{4,2}>0$, while the other two sums
are negative. Therefore, identically to the proof of lower bound in the
previous section, we only need to show that for sufficiently large $K$,
still using $S_{K}=\{n:n^{2}|t-s|\geq K\}$,
\begin{align}
& \sum_{n\in S_{K}}q_{n}(\tilde{I}_{1}+\bar{I}_{1})>2\sum_{n\in
S_{K}}q_{n} \left\vert \tilde{I}_{4,3}+\bar{I}_{4,3} \right\vert ,
\label{I1ge2I43tildes} \\
& \sum_{n\in S_{K}}q_{n}(\tilde{I}_{1}+\bar{I}_{1})>4\sum_{n\in
S_{K}}q_{n} \left\vert \tilde{I}_{4,4}+\bar{I}_{4,4} \right\vert .
\label{I1ge4I44tildes}
\end{align}
This is not difficult. Indeed, we have that both $f$ and $h$ are decreasing,
and for all $u\in \lbrack r,s]$,
\begin{eqnarray*}
\left\vert h\left( u\right) -h\left( r\right) \right\vert &=&\left(
e^{-n^{2}\left( s-u\right) }-e^{-n^{2}\left( s-r\right) }\right) \left\vert
h_{s,t,x,y}\right\vert \\
&\leq &\left( e^{-n^{2}\left( s-u\right) }-e^{-n^{2}\left( s-r\right)
}\right) \left\vert h_{s,t,x,x}\right\vert =\left\vert f\left( u\right)
-f\left( r\right) \right\vert .
\end{eqnarray*}
Hence, exploiting the fact that all the terms in the products defining the $
I_{4,3}$ as well as $\tilde{I}_{4,3}+\bar{I}_{4,3}$ have constant signs, we
can write
\begin{eqnarray*}
\left\vert \tilde{I}_{4,3}+\bar{I}_{4,3} \right\vert
&=&\int_{0}^{s}dr\int_{r}^{s}du\left\vert h(u)-h(r)\right\vert \left\vert 
\frac{\partial K}{\partial u}(u,r)\right\vert
\int_{s}^{t}dv(g(v)-g(r))\left\vert \frac{\partial K}{\partial v}
(v,r)\right\vert \\
&\leq &\int_{0}^{s}dr\int_{r}^{s}du\left\vert f(u)-f(r)\right\vert
\left\vert \frac{\partial K}{\partial u}(u,r)\right\vert
\int_{s}^{t}dv(g(v)-g(r))\left\vert \frac{\partial K}{\partial v}
(v,r)\right\vert \\
&=&\left\vert I_{4,3}\right\vert ,
\end{eqnarray*}
and similarly we get $\left\vert \tilde{I}_{4,4}+\bar{I}_{4,4} \right\vert
\leq \left\vert I_{4,4}\right\vert $. Since the lower bound on $\tilde{I}
_{1}+\bar{I}_{1}$ in (\ref{IIlb}) is as large as the lower bound (\ref
{i1lower}) on $I_{1}$, the proof of the lower bound in the previous section
implies both (\ref{I1ge2I43tildes}) and (\ref{I1ge4I44tildes}), which
finishes the proof of \eqref{claim1} when $h_{s,t,x,y}<0$.\bigskip

\emph{Case }$h_{s,t,x,y}\geq 0$. Here we cannot rely on previous
calculations. Indeed, in this case, $\tilde{I}_{4,3}+\bar{I}_{4,3}\geq 0$
and $\tilde{I}_{4,4}+\bar{I}_{4,4}\geq 0$, while $\tilde{I}_{4,1}+\bar{I}
_{4,1}$ and $\tilde{I}_{4,2}+\bar{I}_{4,2}$ are negative, and we must
therefore control their absolute values. As in the previous case, we only
need to prove that for $K$ large enough, 
\begin{align}
& \sum_{n\in S_{K}}q_{n}(\tilde{I}_{1}+\bar{I}_{1})>2\sum_{n\in
S_{K}}q_{n}\left\vert \tilde{I}_{4,1}+\bar{I}_{4,1}\right\vert ,
\label{I1ge2I41tildes} \\
& \sum_{n\in S_{K}}q_{n}(\tilde{I}_{1}+\bar{I}_{1})>4\sum_{n\in
S_{K}}q_{n}\left\vert \tilde{I}_{4,2}+\bar{I}_{4,2}\right\vert .
\label{I1ge4I42tildes}
\end{align}

Unlike the last section where the full sum had to be invoked to obtain the required lower bounds, here it is possible to prove that the above two inequalities hold without the sums, i.e. for any fixed $n\in S_{K}$. These
fact are established in Appendix \ref{fcc}. 

This proves \eqref{claim1} when $H<\frac12$. The proof of the lemma is thus complete. 
\end{proof}

\begin{prop}
\label{bivariate} Assume hypothesis \textnormal{(\ref{hq})}. Then for all
$t_{0},T,M>0 $, there exists a finite constant $c_{H}>0$ 
such that for all $s,t\in \lbrack t_{0},T]$, $x,y\in S^{1}$ and $
z_{1},z_{2}\in \lbrack -M\,,M]^{d}$, 
\begin{equation*}
p_{t,x;s,y}(z_{1}\,,z_{2})\leq c_{H}\,(\mathbf{\Delta }((t\,,x)\,;(s
\,,y)))^{-d/2}\exp \left( -\frac{\Vert z_{1}-z_{2}\Vert ^{2}}{c_{H}\mathbf{
\Delta }((t\,,x)\,;(s\,,y))}\right) .
\end{equation*}
\end{prop}

\begin{proof}
Let $p_{t,x;s,y}^{i}(\cdot \,,\cdot )$ denote the bivariate density of the
random vector $(u_{i}(t\,,x)\,,u_{i}(s\,,y))$. Note that $
p_{t,x;s,y}^{i}(\cdot \,,\cdot )$ does not depend on $i$.

We follow \cite{Dalang:04} and \cite{Dalang:05}. As in \cite[(3.8)]
{Dalang:04} and \cite[(4.10)]{Dalang:05}, we have that 
\begin{equation}
\begin{split}
p_{t,x;s,y}^{i}(z_{1},z_{2})& \leq \frac{1}{2\pi \sigma _{s,y}\tau }\exp
\left( -\frac{|z_{1}-z_{2}|^{2}}{4\tau ^{2}}\right) \\
& \hskip1.4in\times \exp \left( \frac{|z_{2}|^{2}\,|1-m|^{2}}{4\tau ^{2}}
\right) \exp \left( -\frac{|z_{2}|^{2}}{2\sigma _{s,y}^{2}}\right) ,
\end{split}
\label{biva}
\end{equation}
where 
\begin{equation*}
\begin{split}
& \tau ^{2}:=\sigma _{t,x}^{2}\left( 1-\rho _{t,x;s,y}^{2}\right) ,\qquad
\rho _{t,x;s,y}=\frac{\sigma _{t,x;s,y}}{\sigma _{t,x}\sigma _{s,y}},\qquad
\sigma _{t,x}^{2}={\mathrm{E}}[(u_{i}(t,x))^{2}] \\
& m:=\frac{\sigma _{t,x;s,y}}{\sigma _{s,y}^{2}},\qquad \sigma _{t,x;s,y}=
\text{Cov}\left( u_{i}(t\,,x)\,,u_{i}(s\,,y)\right) .
\end{split}
\end{equation*}

We now show the analogues of (4.12) and Lemma 4.3 in \cite{Dalang:05} in the
case of the fractional heat equation. Fix $s,t\in \lbrack t_{0},T]$, $x,y\in
S^{1}$. We claim that the following estimates hold: 
\begin{align}
|\sigma _{t,x}-\sigma _{s,y}|& \leq c_{H}|t-s|^{2\alpha }.  \label{incr2} \\
c_{H}^{-1}\mathbf{\Delta }((t\,,x)\,;(s\,,y))\leq \sigma _{t,x}^{2}\sigma
_{s,y}^{2}-\sigma _{t,x;s,y}^{2}& \leq c_{H}\,\mathbf{\Delta }
((t\,,x)\,;(s\,,y)),  \label{bisa1} \\
|\sigma _{t,x}^{2}-\sigma _{t,x;s,y}|& \leq c_{H}\left[ \mathbf{\Delta }
((t\,,x)\,;(s\,,y))\right] ^{1/2}.  \label{bisa2}
\end{align}

Indeed, in the proof of Proposition \ref{propt} we have proved that 
\begin{equation*}
{\mathrm{E}}\biggl[\biggl(\int_{0}^{t}e^{-n^{2}\left( t-s\right) }\beta
_{n}^{H}\left( ds\right) \biggr)^{2}\biggr]\leq c_{H}|t-s|^{2\alpha }.
\end{equation*}
Therefore, using \cite[(4.31)]{Dalang:05}, we have 
\begin{equation*}
|\sigma _{t,x}-\sigma _{s,y}|\leq c_{H}\,|\sigma _{t,x}^{2}-\sigma
_{s,y}^{2}|\leq c_{H}|t-s|^{2\alpha }
\end{equation*}%
where $c_{H}$ does not depend on $t\in \lbrack t_{0}\,,T]$. This proves (\ref
{incr2}).

We now prove (\ref{bisa1}). Let $\gamma _{t,x;s,y}^{2}:={\mathrm{E}}
[(u_{i}(t\,,x)-u_{i}(s\,,y))^{2}]$. Then using \cite[(4.42)]{Dalang:05}, 
\begin{equation}
\sigma _{t,x}^{2}\sigma _{s,y}^{2}-\sigma _{t,x;s,y}^{2}=\frac{1}{4}\left(
\gamma _{t,x;s,y}^{2}-(\sigma _{t,x}-\sigma _{s,y})^{2}\right) \left(
(\sigma _{t,x}+\sigma _{s,y})^{2}-\gamma _{t,x;s,y}^{2}\right) .
\label{equat2}
\end{equation}
By Lemma \ref{lemmatex}, $\gamma _{t,x,s,y}^{2}\leq c\mathbf{\Delta }
((t\,,x)\,;(s\,,y))$. Therefore, the second factor of (\ref{equat2}) is
bounded below by a positive constant when $(t\,,x)$ is near $(s\,,y)$.
Furthermore, Lemma \ref{lemmatex} and (\ref{incr2}) yield
\begin{equation*}
\gamma _{t,x,s,y}^{2}-(\sigma _{t,x}-\sigma _{s,y})^{2}\geq c_{H}\mathbf{
\Delta }((t\,,x)\,;(s\,,y)).
\end{equation*}
This proves the lower bound of (\ref{bisa1}) provided $(t\,,x)$ is
sufficiently near $(s\,,y)$.

In order to extend this inequality to all $(t,x)$ and $(s,y)$ in $
[t_{0},T]\times S^{1}$, note that by the contuinuity of the function $
(t,x,s,y)\mapsto \sigma _{t,x}^{2}\sigma _{s,y}^{2}-\sigma _{t,x;s,y}^{2}$,
it suffices to show that 
\begin{equation*}
\sigma _{t,x}^{2}\sigma _{s,y}^{2}-\sigma _{t,x;s,y}^{2}>0\qquad \mbox{if }
(t\,,x)\neq (s\,,y).
\end{equation*}
If this last function was equal to zero there would be $\lambda \in \mathbb{R
}$ such that $u_{i}(t\,,x)=\lambda u_{i}(s\,,y)$ a.s., which is a
contradiction to the lower bound in (\ref{incr}) and the fact that $\mathbf{
\Delta }((t\,,x)\,;(s\,,y))$ is zero only if $(t\,,x)\,=(s\,,y)$. This
completes the proof of the lower bound of (\ref{bisa1}).

In order to prove the upper bound of (\ref{bisa1}), use Lemma \ref{lemmatex} to see
that the first factor in (\ref{equat2}) is bounded above by $c_{H}\mathbf{
\Delta }((t\,,x)\,;(s\,,y))$. As the second factor in (\ref{equat2}) is
bounded above by a constant $c_{H}$, the desired upper bound follows.

It remains to prove (\ref{bisa2}). Use \cite[(4.47)]{Dalang:05} to find 
\begin{equation*}
\begin{split}
|\sigma _{t,x}^{2}-\sigma _{t,x;s,y}|& =\left\vert \gamma _{t,x;s,y}^{2}+
\text{Cov}\left( u_{i}(t\,,x)-u_{i}(s\,,y)\,,u_{i}(s\,,y)\right) \right\vert
\\
& \leq \gamma _{t,x;s,y}^{2}+\gamma _{t,x;s,y}\sigma _{s,y} \\
& \leq c_{H}\,[\mathbf{\Delta }((t\,,x)\,;(s\,,y))]^{1/2},
\end{split}
\end{equation*}
where we have used Lemma \ref{lemmatex} twice in the last inequality. This implies
the desired bound.

Finally, introducing inequalities (\ref{bisa1}) and (\ref
{bisa2}) into (\ref{biva}) and using the independence of the components $
u_{1},...,u_{d}$, the proposition follows.
\end{proof}

\section{Proof of Theorem \protect\ref{t1} and Corollary \protect\ref{c1}}

In order to prove Theorem \ref{t1} we will follow the approach developped in 
\cite{Dalang:05} extended to our situation. For this we shall state and
prove the versions of Theorem 2.1(1), Lemma 2.2(1), Theorem 3.1(1) and Lemma
4.5 in \cite{Dalang:05} needed in our situation.

The first result is an extension of \cite[Lemma 2.2(1)]{Dalang:05} (take $
\alpha=1/2$, $H=1/2$ and $d=\beta$).
\begin{lem}
\label{prel} Let $I$ and $J$ two intervals as in Theorem \ref{t1}. Then for all $N>0$, there exists a finite and positive constant $C=C(I,J,N)$ such that for all 
$a\in[0\,,N]$, 
\begin{equation}  \label{int}
\int_{I} dt \int_{I} ds \int_{J} d x \int_{J} dy \ \frac{e^{-a^2/\mathbf{
\Delta}((t,x);(s,y))}}{\mathbf{\Delta}^{d/2} ((t,x);(s,y))} \le C \, \mathrm{
K}_{d-(\frac{1}{\alpha}+\frac{2}{\alpha \wedge (2H)})}(a),
\end{equation}
where $\mathbf{\Delta}((t,x) \, ; (s,y))$ is the metric defined in \textnormal{(\ref{del})}.
\end{lem} 

\begin{proof}
Write $\alpha_1:= 2 \alpha$ and $\alpha_2:=\alpha \wedge (2H)$. Using the change of variables $\tilde{u}=t-s$ ($t$ fixed), $\tilde{v}=x-y$ ($
x$ fixed) we have that the integral in (\ref{int}) is bounded above by 
\begin{equation*}
4 \vert I \vert \, \vert J \vert \int_0^{\vert I \vert} d\tilde{u}
\int_0^{\vert J \vert} d\tilde{v} \, (\tilde{u}^{\alpha_1}+\tilde{v}
^{\alpha_2})^{-d/2} \exp \biggl(- \frac{a^2}{\tilde{u}^{\alpha_1}+ \tilde{v}
^{\alpha_2}} \biggr).
\end{equation*}
A change of variables [$\tilde{u}^{\alpha_1}=a^2 u$, $\tilde{v}
^{\alpha_2}=a^2 v$] implies that this is equal to 
\begin{equation}  \label{chv}
C a^{\frac{2}{\alpha_1}+\frac{2}{\alpha_2}-d} \int_0^{\vert I \vert^{\alpha_1} a^{-2}} du 
\int_0^{\vert J \vert^{\alpha_2} a^{-2} }
dv\, \frac{u^{\frac{1}{\alpha_1}-1} v^{\frac{1}{\alpha_2}-1}}{(u+v)^{d/2}
} \exp \biggl(-\frac{1}{u+v} \biggr).
\end{equation}
Observe that the last integral is bounded above by 
\begin{equation*}
\int_0^{\vert I \vert^{\alpha_1} a^{-2} } du \int_0^{\vert J \vert^{\alpha_2} a^{-2} } dv \, (u+v )^{\frac{1}{\alpha_1}+\frac{1}{\alpha_2}
-2-\frac{d}{2}} \exp \biggl(-\frac{1}{u+v} \biggr).
\end{equation*}
Pass to polar coordinates to deduce that the preceding is bounded above by $
I_1+I_2(a)$, where 
\begin{equation*}
\begin{split}
&I_1:= \int_0^{K N^{-2}} d\rho \, \rho^{\frac{1}{\alpha_1}+\frac{1
}{\alpha_2}-1-\frac{d}{2}} \exp (-c/\rho), \\
&I_2(a):=\int_{KN^{-2}}^{K a^{-2}} d\rho \, \rho^{\frac{
1}{\alpha_1}+\frac{1}{\alpha_2}-1-\frac{d}{2}},
\end{split}
\end{equation*}
where $K= \vert I \vert^{\alpha_1} \vee \vert J \vert^{\alpha_2}$.
Clearly, $I_1 \leq C < \infty$, and if $d \neq \frac{2}{\alpha_1}+\frac{2}{
\alpha_2}$, then 
\begin{equation*}
I_2(a)=K^{\frac{1}{\alpha_1}+\frac{1}{
\alpha_2}-\frac{d}{2}}  \frac{a^{d-\frac{2}{\alpha_1}-\frac{2}{
\alpha_2}}- N^{ d-\frac{2}{\alpha_1
}-\frac{2}{\alpha_2}}}{\frac{1}{\alpha_1}+\frac{1}{\alpha_2}-
\frac{d}{2}}.
\end{equation*}
If $d > \frac{2}{\alpha_1}+\frac{2}{\alpha_2}$, then $I_2(a) \leq C$ for
all $a \in [0,N]$. If $d< \frac{2}{\alpha_1}+\frac{2}{
\alpha_2}$, then $I_2(a) \leq C a^{d-(\frac{2}{\alpha_1}+\frac{2}{\alpha_2}
)}$. Finally, if $d=\frac{2}{\alpha_1}+
\frac{2}{\alpha_2}$, then 
\begin{equation*}
I_2(a) =2 \left[\ln \left(\frac{1}{a} \right) +\ln (N)\right].
\end{equation*}
Hence, we deduce that for all $a \in [0\,,N]$,
the expression in (\ref{chv}) is bounded above by 
$C \, \mathrm{K}_{d-(\frac{2}{\alpha_1}+\frac{2}{\alpha_2})}(a)$,
provided that $N_0$ in (\ref{k}) is sufficiently large. This proves the lemma.
\end{proof}

The next result uses the proof of \cite[Theorem 2.1(1)]{Dalang:05} applied
to our situation and establishes the lower bound of Theorem \ref{t1}.
\begin{thm}
\label{t1a} Assume hypothesis \textnormal{(\ref{hq})}. 
Let $I\subset (0,T]$ and $J\subset \lbrack 0,2\pi )$ be two fixed
non-trivial compact intervals. Then for all $T>0$ and $M>0$, there exists a finite constant $c_{H}>0$ 
such that for all compact sets $A\subseteq \lbrack -M,M]^{d}$, 
\begin{equation*}
c_{H}\,\textnormal{Cap}_{d-\beta }(A)\leq \P\{u(I\times J)\cap A\neq
\emptyset \},
\end{equation*}
where $\beta :=\frac{1}{\alpha }+(\frac{2}{\alpha }\vee \frac{1}{H})$.
\end{thm}

\begin{proof}
The proof of this result follows exactly the same lines as the proof of \cite
[Theorem 2.1(1)]{Dalang:05}, therefore we will only sketch the steps that
differ. It suffices to replace their $\beta-6$ by our $d-\beta$ with $\beta:=
\frac{1}{\alpha}+(\frac{2}{\alpha} \vee \frac1H)$. Moreover, if $p_{t,x}(y)$
denotes the density of $u(t,x)$ solution of (\ref{equa1}), then we have that
for all $y \in [-M, M]^d$ and $(t,x) \in I \times J$, 
\begin{equation}  \label{strictpos}
p_{t,x}(y)=(2 \pi \sigma^2_{t,x})^{-d/2} e^{-\Vert y \Vert^2/(2
\sigma^2_{t,x})} \geq c_H,
\end{equation}
which proves hypothesis \textbf{A1} of \cite[Theorem 2.1(1)]{Dalang:05}. On
the other hand, our Proposition \ref{bivariate} proves hypothesis \textbf{A2}
with $\mathbf{\Delta}((t,x) \, ; (s,y))$ defined as in (\ref{del}).

We then follow the proof of \cite[Theorem 2.1(1)]{Dalang:05}. 
 Define, for all $z \in \mathbb{R}^d$ and $\epsilon>0$, 
    $\tilde{B}(z \,,\e):=\{y\in\R^d:\ |y-z|<\e\}$, where $|z|:=\max_{1\le j\le d}|z_j|$, and
    \begin{equation*}
        J_\e(z) = \frac{1}{(2\e)^d} \int_{I} dt \, \int_{J } dx \, \1_{\tilde{B}(z,\e)}
        (u(t\,,x)).
    \end{equation*}
In the case $d<\beta$, instead of \cite[(2.31)]{Dalang:05} we will find, using
Proposition \ref{bivariate}, Lemma \ref{prel} and \cite[Lemma 2.3]{Dalang:05}, that for all
$z \in A\subseteq[-M\,,M]^d$ and $\epsilon>0$,
\begin{equation*}
\begin{split}
{\mathrm{E}} \left[ (J_\e(z))^2 \right] \leq &c_H \int_0^{ \vert I \vert} du
\, \int_0^{\vert J \vert} dv \,( u^{2\alpha}+ v^{\alpha \wedge (2H)})^{-d/2} \\
&\leq c_H \int_0^{ \vert I \vert} du \, \Psi_{\vert J \vert, d (\frac{\alpha
}{2} \wedge H)} (u^{d \alpha}) \\
&\leq c_H \int_0^{ \vert I \vert} du \, \mathrm{K}_{1-(\frac{2}{\alpha} \vee
\frac1H)/d}(u^{d \alpha}).
\end{split}
\end{equation*}
We will then consider the different cases: $d< \frac{2}{\alpha} \vee \frac1H$
, $\frac{2}{\alpha} \vee \frac1H <d<\frac{1}{\alpha}+(\frac{2}{\alpha} \vee
\frac1H)$ and $d=\frac{2}{\alpha} \vee \frac1H$. This will prove the case $
d<\beta$.

The case $d \geq \beta$ is proved exactly along the same lines as the proof
of \cite[Theorem 2.1(1)]{Dalang:05}, appealing to (\ref{strictpos}),
Proposition \ref{bivariate} and Lemma \ref{prel}.
\end{proof}

The following result is an extension of \cite[Lemma 4.5]{Dalang:05}.
\begin{lem}
\label{Garsia} Assume hypothesis \textnormal{(\ref{hq})}. For all $p \geq 1 $
, there exists $C_{p,H}>0$ such that for all $\epsilon>0$ and all $(t\,,x)$
fixed, 
\begin{equation}  \label{incr4}
{\mathrm{E}} \left[ \sup_{[\mathbf{\Delta}((t,x) \, ; (s,y))]^{1/2} \leq
\epsilon} \Vert u(t\,,x)-u(s\,,y) \Vert^p \right] \leq C_{p,H} \epsilon^p,
\end{equation}
where $\mathbf{\Delta}((t,x) \, ; (s,y))$ is defined as in \textnormal{(\ref
{del})}.
\end{lem}

\begin{proof}
It suffices to prove (\ref{incr4}) for each coordinate $u_i$, $i=1,\ldots,d$
. We proceed as in \cite[Lemma 4.5]{Dalang:05}, that is, we will use \cite[
Proposition A.1]{Dalang:05} with $S:=S_{\epsilon}=\{(s,y): [\mathbf{\Delta}
((t\,,x)\,;(s\,,y))]^ {1/2} < \epsilon \}$, $\rho((t\,,x)\,,(s\,,y)):=[
\mathbf{\Delta}((t\,,x)\,;(s\,,y))]^{1/2}$, $\mu(dtdx):=dtdx$, $\Psi
(x):=e^{\vert x \vert}-1$, $p(x):=x$, and $f:=u_i$.

Moreover, by (\ref{incr}), the random variable $\mathcal{C}$ defined in \cite
[Proposition A.1]{Dalang:05} satisfies 
\begin{equation*}
{\mathrm{E}} \, [\mathcal{C}]\le {\mathrm{E}} \, \left[ \int_{S_{\epsilon}}
\, dt\, dx \int_{S_{\epsilon}}\, ds\, dy \ \exp\left(\frac{\vert
u_i(t\,,x)-u_i(s\,,y)\vert}{ [\mathbf{\Delta}((t\,,x)\,;(s\,,y))]^{1/2}}
\right)\right] \leq c_H \epsilon^{\beta},
\end{equation*}
where $\beta=\frac{2}{\alpha}+(\frac{4}{\alpha} \vee \frac{2}{H})$.

The rest of the proof follows exactly as in \cite[(4.51)]{Dalang:05} and is
therefore omitted.
\end{proof}

The next result uses the proof of \cite[Theorem 3.1(1)]{Dalang:05} applied
to our situation and establishes the upper bound of Theorem \ref{t1}.
\begin{thm}
\label{t1b} Assume hypothesis \textnormal{(\ref{hq})}. 
Let $I\subset (0,T]$ and $J\subset \lbrack 0,2\pi )$ be two fixed
non-trivial compact intervals. Then for all $T>0$ and $M>0$, there exists a finite constant $c_{H}>0$ 
such that for all Borel sets $A\subseteq \lbrack -M,M]^{d}$, 
\begin{equation*}
\P\{u(I\times J)\cap A\neq \emptyset \}\leq c_{H}\,\mathcal{H}
_{d-\beta }(A),
\end{equation*}
where $\beta :=\frac{1}{\alpha }+(\frac{2}{\alpha }\vee \frac{1}{H})$.
\end{thm}

\begin{proof}
The proof of this result is
similar to the proof of \cite[Theorem 3.1]{Dalang:05}. When $d<\beta $,
there is nothing to prove, so we assume that $d \geq \beta$.

For all positive integers $n$, set $t_k^n:=k 2^{-n/\alpha}$, $
x_{\ell}^n:=\ell 2^{-(2n/\alpha) \vee (n/H)}$, and 
\begin{equation*}
I^n_k = [t_k^n,t_{k+1}^n],\qquad J^n_{\ell} =
[x_{\ell}^n,x_{\ell+1}^n],\qquad R^n_{k,\ell} = I^n_k \times J^n_{\ell}.
\end{equation*}
Then for all $R_{k,\ell}^{n}\subset I\times J$, there exists a constant $c_H>0$
such that the following hitting small balls estimate holds for all
$z\in \mathbb{R}^{d}$ and $\epsilon >0$,
\begin{equation}
\P\{u(R_{k,\ell}^{n})\cap B(z\,,\epsilon )\neq \varnothing \}\leq
c_{H}\epsilon ^{d}.  \label{smallbal}
\end{equation}

Indeed, the proof of (\ref{smallbal}) follows along the same lines
as the proof of \cite[Proposition 4.4]{Dalang:05} for the linear stochastic
heat equation driven by space time white noise.
Namely, consider the random variables
\begin{equation*}\begin{split}
		Y_{k,\ell}^n &:= \inf_{(t,x)\in R_{k,\ell}^n} \left\Vert
			c_{k,\ell}^n(t\,,x) u(t_k^n\,,x^n_{\ell}) - z\right\Vert,\
			\text{and}\\
		Z_{k,\ell}^n &:= \sup_{(t,x)\in R^n_{k,\ell}} \left\Vert
			u(t\,,x) - c^n_{k,\ell}(t\,,x) u(t_k\,,x_{\ell})\right\Vert,
	\end{split}\end{equation*}
where
\begin{equation*}
		c^n_{k,\ell} (t\,,x) := \frac{ \E\left[ u_1(t\,,x)
		u_1(t_k^n\,,x^n_{\ell})\right]}{
		\mathrm{Var}\left[ u_1(t_k^n\,,x^n_{\ell})\right]}.
	\end{equation*}
Note that, because $\{u_i(t,x), u_i(t_k^n, x_{\ell}^n)\}$ is a $2$-dimensional centered Gaussian vector, $Y_{k,\ell}^n$ and $Z_{k,\ell}^n$ are independent.
Hence, the rest of the proof of (\ref{smallbal}) follows as in \cite[Proposition 4.4]{Dalang:05}, using the fact that $\{u_{i}(t,x)\}_{i=1,..,d}$
are independent, centered, Gaussian random variables, with variance bounded
above and below by positive constants, and such
that the upper bound in (\ref{incr}) and Lemma \ref{Garsia} hold.

Now fix $\epsilon \in \,(0\,,1)$ and $n\in \mathbb{N}$ such that $
2^{-n-1}<\epsilon \leq 2^{-n}$, and write 
\begin{equation*}
\P\left\{ u\left( I\times J\right) \cap B(z\,,\epsilon )\neq
\varnothing \right\} \leq \mathop{\sum\sum}_{\substack{ (k,\ell ):  \\ 
R_{k,\ell }^{n}\cap (I\times J)\neq \varnothing }}\P\{u(R_{k,\ell
}^{n})\cap B(z\,,\epsilon )\neq \varnothing \}.
\end{equation*}
The number of pairs $(k\,,\ell )$ involved in the two sums is at most $
2^{\beta n}$, where $\beta :=\frac{1}{\alpha }+(\frac{2}{\alpha }\vee \frac{1
}{H})$. Because $2^{-n-1}<\epsilon $, \eqref{smallbal} implies that 
\begin{equation}
\P\left\{ u\left( I\times J\right) \cap B(z\,,\epsilon )\neq
\varnothing \right\} \leq c_{H}2^{-n(d-\beta )}\leq c_{H}\epsilon ^{d-\beta
},  \label{eq:hit:ball:UB}
\end{equation}
where $c_{H}$ does not depend on $(n\,,\epsilon )$. Therefore, (\ref
{eq:hit:ball:UB}) is valid for all $\epsilon \in \,(0\,,1)$.

Now we use a \emph{covering argument}: Choose $\epsilon \in\,(0\,,1)$ and
let $\{B_i\}_{i=1}^\infty$ be a sequence of open balls in $\mathbb{R}^d$
with respective radii $r_i\in\,(0\,,\epsilon]$ such that 
\begin{equation}  \label{eq:HHHH}
A\subseteq \bigcup_{i=1}^\infty B_i \quad\text{and}\quad \sum_{i=1}^\infty
(2r_i)^{d-\beta} \le \mathcal{H}_{d-\beta} (A)+\epsilon.
\end{equation}

Because $\P\{u(I\times J)\cap A\neq \varnothing \}$ is at most $
\sum_{i=1}^{\infty }\P\{u(I\times J)\cap B_{i}\neq \varnothing \}$, 
\eqref{eq:hit:ball:UB} and \eqref{eq:HHHH} together imply that
\begin{equation*}
\P\left\{ u\left( I\times J\right) \cap A\neq \varnothing \right\}
\leq c_{H}\sum_{i=1}^{\infty }r_{i}^{d-\beta }\leq c_{H}(\mathcal{H}
_{d-\beta }(A)+\epsilon ).
\end{equation*}
Let $\epsilon \rightarrow 0^{+}$ to conclude the proof of the theorem.
\end{proof}

\begin{proof}[Proof of Theorem \protect\ref{t1}]
Theorems \ref{t1a} and \ref{t1b} prove the lower and upper bounds of Theorem 
\ref{t1}, respectively.
\end{proof}

\begin{proof}[Proof of Corollary \protect\ref{c1}]
\begin{itemize}
\item[\textnormal{(a)}] This is an immediate consequence of Theorem \ref{t1}.

\item[\textnormal{(b)}] Let $z \in \mathbb{R}^d$. If $d <\beta$, then $
\textnormal{Cap}_{d-\beta}(\{ z\})=1$. Hence, the lower bound of Theorem \ref
{t1} implies that $\{ z\}$ is not polar. On the other hand, if $d>\beta$,
then $\mathcal{H}_{d-\beta}(\{ z\})=0$ and the upper bound of Theorem \ref
{t1} implies that $\{ z\}$ is polar.

\item[\textnormal{(c)}] Theorem \ref{t1} implies that for $d \geq 1$: $
\textnormal{codim}(u(\mathbb{R}_+ \times\,S^1 ))=(d-\beta)^+$; where $
\textnormal{codim}(E)$ with $E$ a random set is defined in \cite[(5.12)]
{Dalang:05}. Then, when $d>\beta$, \cite[(5.13)]{Dalang:05} implies the
desired result.

The case $d=\beta$ follows using exactly the same argument that lead to the
result in \cite[Corollary 5.3(a)]{Dalang:05} for $d=6$, and is therefore
omitted.
\end{itemize}
\end{proof}

\appendix

\section{Appendix}

\subsection{Riesz-kernel example} \label{a11}

We consider the example of the Riesz kernel. There, we assume that $Q\left(
x\right) =\left\vert x\right\vert ^{-\gamma }$ for some $\gamma \in (0,1)$.
We then first need to show that this is a {\it bonafide} homogeneous spatial
covariance function on the circle (that this is such a function in Euclidean
space is well-known, but here we are restricted to the circle). In other
words, we need to show that
\begin{equation*}
Q\left( x\right) =\sum_{n=0}^{\infty }q_{n}\cos nx,
\end{equation*}
where $\{ q_{n} \}_{n \in \mathbb{N}}$ is a sequence of nonnegative real numbers. Since $Q$ is integrable, we simply
calculate the values $q_{n}$ by (inverse) Fourier transform: using the
symmetry of $Q$, and some scaling, we obtain 
\begin{eqnarray*}
q_{n} &=&\int_{-\pi }^{\pi }e^{inx}\left\vert x\right\vert ^{-\gamma }dx=2\int_{0}^{\pi }\cos \left( nx\right) x^{-\gamma }dx \\
&=&2n^{\gamma -1}\int_{0}^{n\pi }\cos \left( x\right) x^{-\gamma }dx \\
&=&n^{\gamma -1}\sum_{k=0}^{n-1}r\left( k\right),
\end{eqnarray*}
where $r\left( k\right) =2\int_{k\pi }^{\left( k+1\right) \pi }\cos \left(
x\right) x^{-\gamma }dx.$ We can calculate this $r\left( k\right) $ a bit
further: using an integration by parts, we get
\begin{eqnarray*}
r\left( k\right) &=&2\gamma \int_{k\pi }^{(k+1)\pi }x^{-\gamma -1}\sin
\left( x\right) dx \\
&=&2\gamma \left( -1\right) ^{k}\int_{k\pi }^{(k+1)\pi }x^{-\gamma
-1}\left\vert \sin \left( x\right) \right\vert dx.
\end{eqnarray*}
Hence we do indeed have, as announced in the Riesz kernel example, that $
q_{n}=n^{\gamma -1}c\left( n\right) $ where $c\left( n\right) $ is the
partial sum of the alternating sequence with general term $2r\left( k\right) 
$. Also as announced, we clearly see that $r\left( 0\right) >0$, and it is
trivial to prove that $\left\vert r\left( k+1\right) \right\vert <\left\vert
r\left( k\right) \right\vert $, by simply using the change of variable $
x^{\prime }=x-\pi $, and the fact that $\sin \left( x^{\prime }+\pi \right)
=-\sin \left( x^{\prime }\right) $. The partial sums of such an alternating
series are always positive since the first term is positive. All the claims
in the Riesz-kernel example are justified.

\subsection{Fractional Brownian example}\label{a2}

In the fractional noise example, with $H<1/2$ and where $q_{n}=n^{1-2H}$,
the Fourier series representation $Q\left( x\right) =\sum_{n=0}^{\infty
}n^{1-2H}\cos \left( nx\right) $ is only formal because this series diverges
even as an alternating series. Yet we can interpret $B^{H}$ as the spatial
derivative of a space-time fractional Brownian sheet-type process. Indeed,
consider the centered Gaussian field $Y\left( t,x\right) $ which is
fractional Brownian in time with parameter $H$, and has spatial covariance
equal to $R\left( x,y\right) =\left\vert x-y\right\vert ^{2H}$. Using
exactly the same calculations as in the Riesz-kernel case above, but this
time with $\gamma =-2H$, we can still invoke the fact that $x^{-\gamma -1}$
is decreasing, since $2H-1<0$, and thus $R\left( x,y\right) $ can be written
as $\sum_{n=0}^{\infty }c\left( n\right) n^{-2H-1}\cos \left( nx\right) $.
It is then easy to see that $Y$ can be represented as
\begin{equation*}
Y\left( t,x\right) =\sum_{n=0}^{\infty }\sqrt{c\left( n\right) }
n^{-H-1/2}\cos \left( nx\right) B_{n,H}\left( t\right) +\sum_{n=0}^{\infty }
\sqrt{c\left( n\right) }n^{-H-1/2}\sin \left( nx\right) \tilde{B}
_{n,H}\left( t\right)
\end{equation*}
where $\{ B_{n,H} \}_{n \in \mathbb{N}}$ and $\{\tilde{B}_{n,H}\}_{n \in \mathbb{N}}$ are independent sequences of IID
standard fractional Brownian motions. If one then defines the noise in the
heat equation formally (i.e. in the sense of distributions) by
\begin{equation*}
B_{H}\left( t,x\right) =\frac{\partial }{\partial x}Y\left( t,x\right) ,
\end{equation*}
a factor $n$ comes out in the Fourier representation, and one gets that $
B_{H}$ can be written, in the sense of distributions, as
\begin{equation*}
B_{H}\left( t,x\right) =\sum_{n=0}^{\infty }\sqrt{c\left( n\right) }
n^{-H+1/2}\cos \left( nx\right) B_{n,H}\left( t\right) +\sum_{n=0}^{\infty }
\sqrt{c\left( n\right) }n^{-H+1/2}\sin \left( nx\right) \tilde{B}
_{n,H}\left( t\right),
\end{equation*}
from which the formula 
$
q_{n}=c\left( n\right) n^{1-2H}
$
follows, 
i.e. the formal expansion $Q\left( x\right) =\sum_{n=0}^{\infty }c\left(
n\right) n^{-2H+1}\cos \left( nx\right) $ follows immediately. This
justifies using the scale $n^{1-2H}$ to represent the covariance's Fourier
coefficient in this fractional noise case. Note that this justification also
works when $H>1/2$.

It is instructive to note that one can also formally write 
\begin{eqnarray*}
Q\left( x-y\right) &=&\E\left[ \frac{\partial }{\partial x}Y\left(
1,x\right) \frac{\partial }{\partial y}Y\left( 1,y\right) \right] \\
&=&(\partial ^{2}/\partial x\partial y)\left\vert x-y\right\vert
^{2H}=2H\left( 2H-1\right) \left\vert x-y\right\vert ^{2H-2},
\end{eqnarray*}
which is not integrable at the origin ($x=y$) when $H<1/2$, which explains
why one cannot use the pointwise Fourier and/or the Riesz-kernel
representation in this case.

\subsection{Estimates of the kernel $K^{H}$}\label{a3}

We have the following estimates on the kernel $K^H$.
\begin{lem}
\label{lema1} Let $t_0, T \geq 0$ be fixed. Then for any $H< \frac12$ and $
s,t \in [t_0, T]$ with $s \leq t$, there exist positive constants $
c(t_0,T,H) $ and $C(t_0,T,H)$ such that 
\begin{equation*}
\begin{split}
c(t_0,T,H)^{-1} (t-s)^{H-\frac12} &\leq K^H(t,s) \leq c(t_0,T,H)
(t-s)^{H-\frac12} s^{H-\frac12}, \\
C(t_0,T,H)^{-1} (H-\frac12) (t-s)^{H-\frac32} &\leq \frac{\partial{K^H}}{
\partial t} (t,s) \leq C(t_0,T,H) (H-\frac12) (t-s)^{H-\frac32}
\end{split}
\end{equation*}
\end{lem}

\begin{proof}
Theses estimates follow immediately from (\ref{kernel}), (\ref{dkernel}) and 
\cite[Theorem 3.2]{Decreusefond:97}
\end{proof}

The following is a two real variable technical result that is used several
times in this paper.
\begin{lem}
\label{a1} Let $t_0>0$ fixed. Then for any $s\geq t_0$, there exists a
positive constant $c(t_0,H)$ such that 
\begin{equation*}
\int_0^{2n^2s} \biggl(s-\frac{v}{2n^2} \biggr)^{2H-1} v^{2H-1} e^{-v} \, dv
\leq c(t_0,H).
\end{equation*}
\end{lem}

\begin{proof}
We write, following \cite[eq. (25)]{Tindel:03}, 
\begin{equation*}
\begin{split}
& \int_{0}^{2n^{2}s}\biggl(s-\frac{v}{2n^{2}}\biggr)^{2H-1}v^{2H-1}e^{-v}\,dv
\\
& \leq (\frac{s}{2})^{2H-1}\int_{0}^{\infty
}v^{2H-1}e^{-v}\,dv+(n^{2}s)^{2H-1}\int_{n^{2}s}^{2n^{2}s}\biggl(s-\frac{v}{
2n^{2}}\biggr)^{2H-1}e^{-v}\,dv \\
& \leq c_{H}{t_{0}}^{2H-1}+(n^{2}s)^{2H-1}\int_{0}^{n^{2}s}\biggl(\frac{
v^{\prime }}{2n^{2}}\biggr)^{2H-1}e^{-(2n^{2}s-v^{\prime })}\,dv^{\prime } \\
& \leq C(t_{0},H)+c_{H}{t_{0}}^{2H-1}e^{-n^{2}s}(n^{2}s)^{2H} \\
& \leq C(t_{0},H)+C(t_{0},H)\sup_{x\geq s}|e^{-x}x^{2H}| \\
& \leq C(t_{0},H).
\end{split}
\end{equation*}
\end{proof}

\subsection{Further covariance calculations\label{fcc}}

\begin{proof}[Proof of \eqref{I1ge2I41tildes}
]
With the notations of the proof of Lemma \ref{lemmatex}, we will show that for $K$ large enough and for all $n$ such that $n^{2}\left( t-s\right) \geq K$, when $h_{t,s,x,y}\geq 0$,
\begin{equation}
\tilde{I}_{1}+\bar{I}_{1}>2 \left\vert \tilde{I}_{4,1}+\bar{I}
_{4,1} \right\vert .  \label{appendi1}
\end{equation}
This will prove \eqref{I1ge2I41tildes}.

Using Lemma \ref{lema1}, and the trivial bound $h_{t,s,x,y}\leq 2$ applied
to (\ref{hstxy}), we have
\begin{eqnarray*}
\left\vert \tilde{I}_{4,1}+\bar{I}_{4,1} \right\vert
&=&\int_{0}^{s}drK(s,r)h(r)\left( \int_{s}^{t}\left\vert \frac{\partial K}{
\partial u}\left( u,r\right) \right\vert du\right) g(r) \\
&\leq &c_{H}\int_{0}^{s}dr\left( s-r\right) ^{H-1/2}e^{-n^{2}\left(
t+s-2r\right) }\left( \left( s-r\right) ^{H-1/2}-\left( t-r\right)
^{H-1/2}\right) \\
&=&c_{H}e^{-n^{2}\left( t-s\right) }\int_{0}^{s}dr\ r^{H-1/2}\left(
r^{H-1/2}-\left( r+t-s\right) ^{H-1/2}\right) e^{-2n^{2}r}.
\end{eqnarray*}
We evaluate the integral above by splitting it up according to whether $r$
exceeds $n^{-2}$. We also assume that $n^{2}\left( t-s\right) \geq 1$, i.e.
we restrict $K\geq 1$. Hence
\begin{eqnarray*}
&&\int_{0}^{n^{-2}}dr\ r^{H-1/2}\left( r^{H-1/2}-\left( r+t-s\right)
^{H-1/2}\right) e^{-2n^{2}r} \\
&& \qquad \qquad \qquad \leq \int_{0}^{n^{-2}}dr\ r^{H-1/2}\left( r^{H-1/2}-\left( 2t-2s\right)
^{H-1/2}\right) \\
&&\qquad \qquad \qquad  =\int_{0}^{n^{-2}}dr\ \left( r^{2H-1}-r^{H-1/2}\left( 2t-2s\right)
^{H-1/2}\right) \\
&&  \qquad \qquad \qquad \leq c_{H}n^{-4H}.
\end{eqnarray*}
The other piece is
\begin{eqnarray*}
&&\int_{n^{-2}}^{s}dr\ r^{H-1/2}\left( r^{H-1/2}-\left( r+t-s\right)
^{H-1/2}\right) e^{-2n^{2}r} \\
&& \qquad \leq c_{H}\int_{n^{-2}}^{s}dr\ r^{H-1/2}\left( t-s\right)
r^{H-3/2}e^{-2n^{2}r}=c_{H}\left( t-s\right) \int_{n^{-2}}^{s}dr\
r^{2H-2}e^{-2n^{2}r} \\
&& \qquad =c_{H}n^{-2}n^{4-4H}\left( t-s\right) \int_{1}^{n^{2}s}dx\ x^{2H-2}e^{-2x}
\\
&& \qquad \leq c_{H}n^{-4H}n^{2}\left( t-s\right) \int_{1}^{\infty }dx\
x^{2H-2}e^{-2x} \\
&&\qquad  =c_{H}n^{-4H}n^{2}\left( t-s\right) .
\end{eqnarray*}
In conclusion, we get 
\begin{equation*}
\left\vert \tilde{I}_{4,1}+\bar{I}_{4,1} \right\vert \leq c^1_{H}n^{-4H}\left(
1+n^{2}\left( t-s\right) \right) e^{-n^{2}\left( t-s\right) }.
\end{equation*}
Since the function $x\mapsto \left( 1+x\right) e^{-x}$ decreases to $0$ as $
x $ increases to $\infty $, we only need to choose $K$ sufficiently large such that for all $n$ with
$n^{2}\left( t-s\right) \geq K$, $\left\vert \tilde{I}
_{4,1}+\bar{I}_{4,1}\right\vert \leq 2^{-1}c^1_{H}n^{-4H}(
1-e^{-n^2(t-s)}) ^{2}\leq \tilde{I}_{1}+\bar{I}_{1}$, where $c^1_H$ is the constant in \eqref{IIlb}. This completes the
proof of (\ref{appendi1}).
\end{proof}

\begin{proof}[Proof of \eqref{I1ge4I42tildes}]
We now show that for $K$ large enough and for all $n$ such that $n^{2}\left(
t-s\right) \geq K$, when $h_{t,s,x,y}\geq 0$,
\begin{equation}
\tilde{I}_{1}+\bar{I}_{1}>2 \left\vert \tilde{I}_{4,2}+\bar{I}
_{4,2} \right\vert.  \label{appendi2}
\end{equation}
This will prove \eqref{I1ge4I42tildes}.

Again using Lemma \ref{lema1}, and the bound $h_{t,s,x,y}\leq 2$ applied to (\ref{hstxy}), we have
\begin{eqnarray*}
\left\vert \tilde{I}_{4,2}+\bar{I}_{4,2} \right\vert
&=&h_{t,s,x,y}\int_{0}^{s}drK(s,r)e^{-n^{2}\left( s-r\right)
}\int_{s}^{t}du(g(u)-g(r))\left\vert \frac{\partial K}{\partial u}
(u,r)\right\vert \\
&\leq &c_{H}\int_{0}^{s}dr (s-r)^{H-1/2}e^{-n^{2}\left(
s+t-r\right) }\int_{s}^{t}du\left( e^{n^{2}u}-e^{n^{2}r}\right) (u-r)^{H-3/2}.
\end{eqnarray*}
We cut this integral into three pieces. First calculate the piece for $u>s+n^{-2}$:
\begin{eqnarray*}
&&\int_{0}^{s}dr (s-r)^{H-1/2}e^{-n^{2}\left(
s+t-r\right) }\int_{s+n^{-2}}^{t}du\left( e^{n^{2}u}-e^{n^{2}r}\right)
(u-r)^{H-3/2} \\
&& \qquad \leq \int_{0}^{s}dr (s-r)^{H-1/2}e^{-n^{2}\left(
s+t-2r\right) }\int_{s+n^{-2}}^{t}due^{n^{2}(u-r) }(u-r)^{H-3/2} \\
&& \qquad =n^{-2H+1}\int_{0}^{s}dr (s-r)^{H-1/2}e^{-n^{2}\left(
s+t-2r\right) }\int_{\left( s-r\right) n^{2}+1}^{\left( t-r\right)
n^{2}}e^{x}x^{H-3/2}dx \\
&& \qquad =n^{-4H}\int_{0}^{sn^{2}}dy\ y^{H-1/2}e^{-y}e^{-n^{2}\left( t-s\right)
}\int_{y+1}^{y+n^{2}\left( t-s\right) }e^{x}x^{H-3/2}dx.
\end{eqnarray*}
Now, for any fixed constants $y_{0}\left( H\right)$ and $y_{1}\left( H\right)$ such that $y_{1}
>y_{0} +1$, the above term with
the $y$-integral restricted to $y\leq y_{0}$ can be written
as follows:
\begin{eqnarray*}
&&n^{-4H}\int_{0}^{y_{0}}dy\ y^{H-1/2}e^{-y}e^{-n^{2}\left( t-s\right)
}\left( \int_{y+1}^{y_{1}}e^{x}x^{H-3/2}dx+\int_{y_{1}}^{y+n^{2}\left(
t-s\right) }e^{x}x^{H-3/2}dx\right) \\
&& \qquad\leq n^{-4H}\int_{0}^{y_{0}}dy\ y^{H-1/2}\left( e^{-n^{2}\left( t-s\right)
}c\left( H,y_{1}\right) +y_{1}^{H-3/2
}e^{y_{0}}\right).
\end{eqnarray*}
We now choose $y_{1}$ and $K$ large enough such
that for all $n$ with $n^{2}\left( t-s\right) \geq K$ and for any choice of $y_0$, 
the above equation is smaller than $c_{H}n^{-4H}$ with $c_H \leq 2^{-1} c^1_{H} (
1-e^{-n^2(t-s)}) ^{2}$, where $c^1_H$ is the constant in \eqref{IIlb}.

For the other part of the integral in $y$ we get
\begin{eqnarray*}
&&n^{-4H}\int_{y_{0}}^{sn^{2}}dy\ y^{H-1/2}e^{-y}e^{-n^{2}\left( t-s\right)
}\int_{y+1}^{y+n^{2}\left( t-s\right) }e^{x}x^{H-3/2}dx \\
&& \qquad \leq n^{-4H}\int_{y_{0}}^{sn^{2}}dy\ y^{2H-2}e^{-y}e^{-n^{2}\left(
t-s\right) }\int_{y+1}^{y+n^{2}\left( t-s\right) }e^{x}dx \\
&& \qquad \leq n^{-4H}\int_{y_{0}}^{sn^{2}}dy\ y^{2H-2} \\
&& \qquad \leq c_{H}n^{-4H} y_{0}^{2H-1},
\end{eqnarray*}
and it is sufficient to take $y_{0}$ large enough to ensure that this last
expression is smaller than $c_H n^{-4H}$ with  $c_H \leq 2^{-1} c^1_{H} (
1-e^{-n^2(t-s)}) ^{2}$.

Now we calculate the piece for $u\in \lbrack s,s+n^{-2}]$ and $r\in \lbrack
s-n^{-2},s]$. This yields a piece bounded above by
\begin{eqnarray*}
&&c_{H}\int_{s-n^{-2}}^{s}dr (s-r)
^{H-1/2}e^{-n^{2}t}\int_{s}^{s+n^{-2}}du\left(
e^{n^{2}s+1}-e^{n^{2}s-1}\right) (u-r)^{H-3/2} \\
&& \qquad \leq c_{H}e^{-n^{2}\left( t-s\right) }\int_{s-n^{-2}}^{s}dr(s-r)^{H-1/2}\left( (s-r)^{H-1/2}-
(s-r+n^{-2})^{H-1/2}\right) \\
&& \qquad =c_{H}e^{-n^{2}\left( t-s\right) }n^{-4H}\int_{0}^{1}x^{H-1/2}\left(
x^{H-1/2}-\left( x+1\right) ^{H-1/2}\right) dx \\
&& \qquad =c_{H}e^{-n^{2}\left( t-s\right) }n^{-4H}
\end{eqnarray*}
which can obviously be made smaller than $2^{-1} c^1_{H} (
1-e^{-n^2(t-s)}) ^{2}$, for all $n$ such that $
n^{2}\left( t-s\right) \geq K$, provided that $K$ is large enough.

The last piece to deal with is
\begin{eqnarray*}
&&c_{H}\int_{0}^{s-n^{-2}}dr (s-r)
^{H-1/2}e^{-n^{2}\left( s+t-r\right) }\int_{s}^{s+n^{-2}}du\left(
e^{n^{2}u}-e^{n^{2}r}\right) (u-r)^{H-3/2} \\
&& \qquad \leq c_{H}\int_{0}^{s-n^{-2}}dr (s-r)
^{H-1/2}e^{-n^{2}\left( s+t-r\right) }\int_{s}^{s+n^{-2}}du\
e^{n^{2}u}(u-r)^{H-3/2} \\
&& \qquad \leq c_{H}e\int_{0}^{s-n^{-2}}dr (s-r)
^{H-1/2}e^{-n^{2}\left( t-r\right) }\int_{s}^{s+n^{-2}}du\ 
(u-r)^{H-3/2} \\
&& \qquad =c_{H}e^{-n^{2}\left( t-s\right) }\int_{0}^{s-n^{-2}}dr (
s-r)^{H-1/2}\left( (s-r)^{H-1/2}-(
s-r+n^{-2})^{H-1/2}\right) \\
&& \qquad \leq c_{H}e^{-n^{2}\left( t-s\right) }n^{-4H}\int_{1}^{\infty
}x^{H-1/2}\left( x^{H-1/2}-\left( x+1\right) ^{H-1/2}\right) dx \\
&& \qquad \leq c_{H}e^{-n^{2}\left( t-s\right) }n^{-4H}\int_{1}^{\infty }x^{H-3/2}dx
\\
&& \qquad =c_{H}e^{-n^{2}\left( t-s\right) }n^{-4H},
\end{eqnarray*}
and the conclusion is the same as before. This finishes the proof of (\ref
{appendi2}).
\end{proof}

\end{document}